\documentclass[11pt,reqno]{amsart}
\usepackage{amssymb,amscd,amsbsy}
\usepackage{amssymb,amscd,amsbsy,mathrsfs}
\newcommand{\lb}{\linebreak}

\renewcommand{\a}{\alpha}
\renewcommand{\b}{\beta}

\renewcommand{\d}{\delta}
\newcommand{\e}{\varepsilon}

\newcommand{\z}{\zeta}

\renewcommand{\l}{\lambda}

\newcommand{\s}{\sigma}

\newcommand{\f}{\varphi}
\renewcommand{\o}{\omega}

\newcommand{\D}{\Delta}
\renewcommand{\L}{\Lambda}

\renewcommand{\O}{\Omega}

\newcommand{\B}{{\mathscr B}}

\newcommand{\E}{{\mathscr E}}
\newcommand{\cd}{{\mathscr D}}
\newcommand{\F}{{\mathscr F}}

\newcommand{\h}{{\mathscr H}}

\newcommand{\X}{{\mathscr X}}
\newcommand{\Y}{{\mathscr Y}}

\newcommand{\C}{{\Bbb C}}

\newcommand{\R}{{\Bbb R}}
\newcommand{\Z}{{\Bbb Z}}

\newcommand{\0}{{\boldsymbol{0}}}

\newcommand{\bs}{\boldsymbol}

\newcommand{\bS}{{\boldsymbol S}}

\newcommand{\rf}[1]{(\ref{#1})}

\newcommand{\df}{\stackrel{\mathrm{def}}{=}}

\newcommand{\re}{\operatorname{Re}}

\newcommand{\supp}{\operatorname{supp}}

\newcommand{\rank}{\operatorname{rank}}
\newcommand{\const}{\operatorname{const}}

\newcommand{\eeq}{\end{equation}}
\newcommand{\beq}{\begin{equation}}
\newcommand{\bay}{\begin{eqnarray}}
\newcommand{\ba}{\begin{align*}}
\newcommand{\ea}{\end{align*}}
\newcommand{\ey}{\end{eqnarray}}
\newcommand{\bey}{\begin{eqnarray*}}
\newcommand{\eey}{\end{eqnarray*}}

\newcommand{\imp}{\Rightarrow}
\newcommand{\be}{\infty}

\newcommand{\bl}{\blacksquare}

\newcommand{\Pf}{{\bf Proof. }}
\newcommand{\im}{\operatorname{Im}}
\renewcommand{\re}{\operatorname{Re}}
\newcommand{\ov}{\overline}

\newtheorem{thm}{\hspace{\parindent}Theorem}[section]

\newtheorem{cor}[thm]{\hspace{\parindent}Corollary}
\newtheorem{lem}[thm]{\hspace{\parindent}Lemma}

\theoremstyle{remark}

\newtheorem*{rem*}{Remark}

\newcommand\Li{{\rm Lip}}
\newcommand\fM{\frak M}

\newcommand\dg{\frak D}

\newcommand{\fI}{{\frak I}}

\newcommand{\qm}{\quad\mbox{and}\quad}

\begin{document}

\newcommand{\vse}{\vspace{.2in}}
\numberwithin{equation}{section}

\title{Functions of normal operators under perturbations}
\author{A.B. Aleksandrov, V.V. Peller, D.S. Potapov, and F.A. Sukochev}

\begin{abstract}
In \cite{Pe1}, \cite{Pe2}, \cite{AP1}, \cite{AP2}, and \cite{AP3} sharp estimates for $f(A)-f(B)$ were obtained for self-adjoint operators $A$ and $B$ and for various classes of functions $f$ on the real line $\R$.
In this paper we extend those results to the case of functions of normal operators. We show that if a function $f$ belongs to the H\"older class $\L_\a(\R^2)$, $0<\a<1$, of functions of two variables, and $N_1$ and $N_2$ are normal operators, then
$\|f(N_1)-f(N_2)\|\le\const\|f\|_{\L_\a}\|N_1-N_2\|^\a$. We obtain a more general result for functions in the space
$\L_\o(\R^2)=\big\{f:~|f(\z_1)-f(\z_2)|\le\const\o(|\z_1-\z_2|)\big\}$ for an arbitrary modulus of continuity $\o$. We prove that if $f$ belongs to the Besov class $B_{\be1}^1(\R^2)$, then it is operator Lipschitz, i.e.,
$\|f(N_1)-f(N_2)\|\le\const\|f\|_{B_{\be1}^1}\|N_1-N_2\|$. We also study properties of $f(N_1)-f(N_2)$ in the case when $f\in\L_\a(\R^2)$ and $N_1-N_2$ belongs to the Schatten-von Neuman class $\bS_p$.
\end{abstract}

\maketitle

\

\begin{center}
{\Large Contents}
\end{center}

\

\begin{enumerate}
\item[1.] Introduction \quad\dotfill \pageref{In}
\item[2.]  Function spaces \quad\dotfill \pageref{Fsp}
\item[3.] Operator ideals \quad\dotfill \pageref{Opid}
\item[4.] Double operator integrals  \quad\dotfill \pageref{dois}
\item[5.] The basic formula in terms of double operator integrals \quad\dotfill \pageref{bfdoi}
\item[6.] Proof of Theorem \ref{SMdd} \quad\dotfill \pageref{PoT}
\item[7.] Operator Lipschitzness and preservation of operator ideals \quad\dotfill \pageref{OLpoi}
\item[8.]  Operator H\"older functions and arbitrary moduli of continuity \quad\dotfill \pageref{oHfamc}
\item[9.]  Perturbations of class $\bS_p$ and more general operator ideals\quad\dotfill \pageref{pocS}
\item[10.]  Commutators and quasicommutators\quad\dotfill \pageref{Caq}
\item[] References \quad\dotfill \pageref{bibl}
\end{enumerate}

\

\setcounter{section}{0}
\section{\bf Introduction}
\setcounter{equation}{0}
\label{In}

\medskip

The purpose of this paper is to generalize results of the papers \cite{Pe1}, \cite{Pe2}, \cite{AP1}, \cite{AP2}, and \cite{AP3} to the case of normal operators.

A Lipschitz function $f$ on the real line $\R$ (i.e., a function satisfying
the inequality $|f(x)-f(y)|\le\const|x-y|$, $x,\,y\in\R$) does not have to be {\it operator Lipschitz}.  In other words, a Lipschitz function $f$ does not necessarily satisfy the inequality
$$
\|f(A)-f(B)\|\le\const\|A-B\|
$$
for arbitrary self-adjoint operators $A$ and $B$ on Hilbert space.
The existence of such functions was proved in \cite{F1}.
Later Kato proved in \cite{K} that the function $f(x)=|x|$ is not operator Lipschitz. Note also that earlier McIntosh established in \cite{Mc} a similar result for commutators (i.e., the function $f(x)=|x|$ is not commutator Lipschitz).  

In \cite{Pe2} and \cite{Pe3} necessary conditions were found for a function $f$ to be operator Lipschitz. In particular, it was shown in \cite{Pe2} that if $f$ is operator Lipschitz, then $f$ belongs locally to the Besov space $B_{11}^1(\R)$. This also implies that Lipschitz  functions do not have to be operator Lipschitz. Note that in \cite{Pe2} and \cite{Pe3} stronger necessary conditions were also obtained.
Note also that the necessary conditions obtained in \cite{Pe1} and \cite{Pe2} are based on the trace class criterion for Hankel operators, see \cite{Pe1} and \cite{Pe4}, Ch. 6.

On the other hand, it was shown in \cite{Pe2} and \cite{Pe3} that if $f$ belongs to the Besov class $B_{\be1}^1(\R)$, then $f$ is operator Lipschitz. We refer the reader to \cite{Pee} for information on Besov spaces. 

It was shown in \cite{AP1} and \cite{AP2}
that the situation dramatically changes if we consider H\"older classes $\L_\a(\R)$ with $0<\a<1$. In this case such functions are necessarily {\it operator H\"older of order $\a$}, i.e., the condition
$|f(x)-f(y)|\le\const|x-y|^\a$, $x,\,y\in\R$,
implies that for self-adjoint operators
$A$ and $B$ on Hilbert space,
$$
\|f(A)-f(B)\|\le\const\|A-B\|^\a.
$$
Another proof of this result was found in \cite{FN2}.

This result was generalized in \cite{AP1} and \cite{AP2} to the case of functions of class $\L_\o(\R)$ for arbitrary moduli of continuity $\o$. This class consists of functions $f$ on $\R$, for which $|f(x)-f(y)|\le\const\o(|x-y|)$, $x,\,y\in\R$. 

Let us also mention that in \cite{AP1} and \cite{AP3} properties of operators $f(A)-f(B)$
were studied for functions $f$ in $\L_\a(\R)$ and self-adjoint operators $A$ and $B$ whose difference $A-B$ belongs to Schatten--von Neumann classes $\bS_p$.

In \cite{AP1}, \cite{AP2} and \cite{AP4} analogs of the above results were obtained for higher order operator differences. 

We also mention here that the papers \cite{AP1}, \cite{AP2}, \cite{AP3}, \cite{AP4}, \cite{AP5}, and \cite{Pe5} study problems of perturbation theory for unitary operators, contractions, and dissipative operators. 

In this paper we are going to study the case of (not necessarily bounded) normal operators. 

In \S\,\ref{OLpoi} we prove that if $f$ is a function on $\R^2$ that belongs to the Besov class $B_{\be1}^1\big(\R^2\big)$, then it is an {\it operator Lipschitz function on} $\R^2$, i.e., 
$$
\big\|f(N_1)-f(N_2)\big\|\le\const\|N_1-N_2\|
$$
for arbitrary normal operators $N_1$ and $N_2$. Note that we say that the operator $N_1-N_2$ is bounded if the domains $\cd_{N_1}$ and $\cd_{N_2}$ of $N_1$ and $N_2$
coincide and $N_1-N_2$ is bounded on $\cd_{N_1}$. If $N_1-N_2$ is not a bounded operator, we say that $\|N_1-N_2\|=\be$.

Note, however, that the proof of the corresponding result for self-adjoint operators obtained in \cite{Pe3} does not work in the case of normal operators. In the case of self-adjoint operators it was shown in \cite{Pe3} that for functions $f$ in the Besov space $B_{\be1}^1(\R)$ and self-adjoint operators $A$ and $B$ with bounded $A-B$, the following formula holds:
$$
f(A)-f(B)=\iint\limits_{\R\times\R}\frac{f(x)-f(y)}{x-y}\,dE_A(x)(A-B)\,dE_B(y).
$$
The expression on the right is a double operator integral.
However, in the case of normal operators a similar formula holds for arbitrary normal operators only for linear functions (see a more detailed discussion in \S\,5).

In \S\,\ref{bfdoi} we obtain a new formula for $f(N_1)-f(N_2)$ in terms of double operator integrals for suitable functions $f$ on $\C$ and normal operators $N_1$ and $N_2$ with bounded $N_1-N_2$. The validity of this formula depends on the fact that certain divided differences are Schur multipliers. This will be proved in \S\,\ref{PoT}.

In \S\,\ref{oHfamc} we prove that as in the case of self-adjoint operators, 
H\"older functions of order $\a$, $0<\a<1$, must be operator H\"older of order $\a$. We also consider the case of arbitrary moduli of continuity. Note that in \cite{FN1} some weaker results were obtained.

Section \ref{pocS} is devoted to the study of properties of $f(N_1)-f(N_2)$,
where $N_1$ and $N_2$ are normal operators whose difference $N_1-N_2$ belongs to the Schatten--von Neumann class $\bS_p$ and $f$ belongs to the H\"older class $\L_\a\big(\R^2\big)$. We obtain analogs for normal operators of the results of \cite{AP1} and \cite{AP2} for self-adjoint operators. We also obtain much more general results for normal operators $N_1$ and $N_2$ whose difference $N_1-N_2$ belongs to ideals of operators on Hilbert space.

Finally, in \S\,\ref{Caq} we obtain estimates for quasicommutators $f(N_1)R-Rf(N_2)$ in terms of $N_1R-RN_2$ and $N_1^*R-RN_2^*$. 

In \S\,\ref{Fsp} we give a brief introduction to Besov spaces and the spaces $\L_\o\big(\R^2\big)$ of functions of two real variables. In \S\,\ref{Opid} we review ideals of operators on Hilbert space. Finally, \S\,\ref{dois} is an
introduction to the Birman--Solomyak theory of double operator integrals.

Note that the results of this paper were announced in the note \cite{APPS}.

Throughout the paper we identify the complex plane 
$\C$ with $\R^2$.

\

\section{\bf Function spaces}
\setcounter{equation}{0}
\label{Fsp}

\

In this section we collect necessary information on Besov spaces and the spaces $\L_\o\big(\R^2\big)$ of functions of two real variables.

\medskip

{\bf 2.1. Besov classes.}
The purpose of this subsection is to give a brief introduction to Besov spaces that play an important role in problems of perturbation theory.
We need the Besov spaces on $\R^2$ only.

%Let $1\le p,\,q\le\be$ and $s\in\R$. The Besov class $B^s_{pq}$ of functions (or
%distributions) on $\T$ can be defined in the following way.

Let $w$ be an infinitely differentiable function on $\R$ such
that
\bay
\label{w}
w\ge0,\quad\supp w\subset\left[\frac12,2\right],\quad\mbox{and} \quad w(x)=1-w\left(\frac x2\right)\quad\mbox{for}\quad x\in[1,2].
\ey

We define the functions $W_n$ on $\R^2$ by
$$
\F W_n(x)=w\left(\frac{|x|}{2^n}\right),\quad n\in\Z, \quad x=(x_1,x_2),
\quad|x|\df\big(x_1^2+x_2^2\big)^{1/2},
$$
where $\F$ is the {\it Fourier transform} defined on $L^1\big(\R^2\big)$ by
$$
\big(\F f\big)(t)=\int_{\R^2} f(x)e^{-{\rm i}(x,t)}\,dx,\quad x=(x_1,x_2),
\quad t=(t_1,t_2), \quad(x,t)\df x_1t_1+x_2t_2.
$$

With each tempered distribution $f\in{\mathscr S}^\prime\big(\R^2\big)$, we
associate a sequence $\{f_n\}_{n\in\Z}$,
\bay
\label{fn}
f_n\df f*W_n.
\ey
Initially we define the (homogeneous) Besov class $\dot B^s_{pq}\big(\R^2\big)$,
$s>0$, $1\le p,\,q\le\be$, as the space of all
$f\in{\mathscr S}^\prime(\R^2)$
such that
\bay
\label{Wn}
\{2^{ns}\|f_n\|_{L^p}\}_{n\in\Z}\in\ell^q(\Z).
\ey
According to this definition, the space $\dot B^s_{pq}(\R^2)$ contains all polynomials. Moreover, the distribution $f$ is defined by the sequence $\{f_n\}_{n\in\Z}$
uniquely up to a polynomial. It is easy to see that the series $\sum_{n\ge0}f_n$ converges in ${\mathscr S}^\prime(\R)$.
However, the series $\sum_{n<0}f_n$ can diverge in general. It is easy to prove that the
series 
\bay
\label{ryad}
\sum_{n<0}\frac{\partial^r f_n}{\partial x_1^k\partial x_2^{r-k}}
\ey
converges uniformly on $\R^2$ for every nonnegative integer
$r>s-2/p$ and $0\le k\le r$. Note that in the case $q=1$ the series \rf{ryad}
converges uniformly, whenever $r\ge s-2/p$ and $0\le k\le r$.

Now we can define the modified (homogeneous) Besov class $B^s_{pq}\big(\R^2\big)$. We say that a distribution $f$
belongs to $B^s_{pq}(\R^2)$ if \rf{Wn} holds and
$$
\frac{\partial^r f}{\partial x_1^k\partial x_2^{r-k}}=\sum_{n\in\Z}\frac{\partial^r f_n}{\partial x_1^k\partial x_2^{r-k}}
$$
in the space ${\mathscr S}^\prime\big(\R^2\big)$, where $r$ is
the minimal nonnegative integer such that $r>s-2/p$ ($r\ge s-2/p$ if $q=1$) and $0\le k\le r$. Now the function $f$ is determined uniquely by the sequence $\{f_n\}_{n\in\Z}$ up
to a polynomial of degree less than $r$, and a polynomial $\varphi$ belongs to $B^s_{pq}\big(\R^2\big)$
if and only if $\deg\varphi<r$.

To define a regularized de la Vall\'ee Poussin type kernel $V_n$, we define the $C^\be$ function $v$ on $\R$ by
\bey
%\label{VP}
v(x)=1\quad\mbox{for}\quad x\in[-1,1]\quad\mbox{and}\quad v(x)=w(|x|)\quad\mbox{if}\quad |x|\ge1,
\eey
where $w$ is the function defined by \rf{w}. Now
we can define the de la Vall\'ee Poussin type functions $V_n$ by
$$
\F V_n(x)=v\left(\frac{|x|}{2^n}\right),\quad n\in\Z,\quad x=(x_1,x_2).
$$
We put $V\df V_0$. Clearly, $V_n(x)=2^{2n}V(2^n x)$.

Besov classes admit many other descriptions.
We give here the definition in terms of finite differences.
For $h\in\R^2$, we define the difference operator $\D_h$,
$$
(\D_hf)(x)=f(x+h)-f(x),\quad x\in\R^2.
$$
It is easy to see that $B_{pq}^s\big(\R^2\big)\subset L^1_{\rm loc}\big(\R^2\big)$ for every $s>0$
and $B_{pq}^s\big(\R^2\big)\subset C\big(\R^2\big)$ for every $s>2/p$. Let $s>0$ and let $m$ be a positive integer such that $m-1\le s<m$. 
The Besov space $B_{pq}^s\big(\R^2\big)$ can be defined as the set of
functions $f\in L^1_{\rm loc}\big(\R^2\big)$ such that
$$
\int_{\R^2}|h|^{-2-sq}\|\D^m_h f\|_{L^p}^q\,dh<\be\quad\mbox{for}\quad q<\be
$$
and
\bay
\label{pbe}
\sup_{h\not=0}\frac{\|\D^m_h f\|_{L^p}}{|h|^s}<\be\quad\mbox{for}\quad q=\be.
\ey
However, with this definition the Besov space can contain polynomials of higher degree than in the case of the first definition given above. 

We use the notation $B_p^s\big(\R^2\big)$ for $B_{pp}^s\big(\R^2\big)$.

For $\a>0$, denote by $\L_\a\big(\R^2\big)$ the H\"older--Zygmund class that consists
of functions $f\in C\big(\R^2\big)$ such that 
$$
|(\D_h^mf)(x)|\le\const|h|^\a,\quad x,~h\in\R^2,
$$
where $m$ is the smallest integer greater than $\a$.
By \rf{pbe}, we have $\L_\a\big(\R^2\big)=B^\a_\be\big(\R^2\big)$.

We refer the reader to \cite{Pee} and \cite{T} for more detailed information on Besov spaces.

\medskip

{\bf 2.2. Spaces $\bs{\L_\o}\big(\R^2\big)$.} Let $\o$ be a modulus of continuity, i.e., $\o$ is a nondecreasing continuous function on $[0,\be)$
such that $\o(0)=0$, $\o(x)>0$ for $x>0$, and
$$
\o(x+y)\le\o(x)+\o(y),\quad x,~y\in[0,\be).
$$
We denote by $\L_\o\big(\R^2\big)$ the space of functions on $\R^2$ such that
$$
\|f\|_{\L_\o(\R^2)}\df\sup_{x\ne y}\frac{|f(x)-f(y)|}{\o(|x-y|)}<\be.
$$

\begin{thm}
\label{Vn}
There exists a constant $c>0$ such that for an arbitrary
modulus of continuity $\o$ and for an arbitrary function $f$ in $\L_\o\big(\R^2\big)$,
the following inequality holds:
\bay
\label{VPn}
\|f-f*V_n\|_{L^\be}\le c\,\o\big(2^{-n}\big)\|f\|_{\L_\o(\R^2)},\quad n\in\Z.
\ey
\end{thm}

\Pf We have
\begin{align*}
\big|f(x)-\big(f*V_n\big)(x)\big|&=2^{2n}\left|\int_{\R^2}\big(f(x)-f(x-y)\big)V\left(2^ny\right)\,dy\right|\\[.2cm]
&\le2^{2n}\|f\|_{\L_\o(\R^2)}\int_{\R^2}\o(|y|)\,\left|V\left(2^ny\right)\right|\,dy\\[.2cm]
&=2^{2n}\|f\|_{\L_\o(\R^2)}\int_{\{|y|\le 2^{-n}\}}\o(|y|)\,\left|V\left(2^ny\right)\right|\,dy\\[.2cm]
&+2^{2n}\|f\|_{\L_\o(\R^2)}\int_{\{|y|>2^{-n}\}}\o(|y|)\,\left|V\left(2^ny\right)\right|\,dy.
\end{align*}
Clearly,
$$
2^{2n}\int_{\{|y|\le 2^{-n}\}}\o(|y|)\,\left|V\left(2^ny\right)\right|\,dy
\le\o\big(2^{-n}\big)\|V\|_{L^1}.
$$
On the other hand, keeping in mind the obvious inequality
$2^{-n}\o(|y|)\le2|y|\o\big(2^{-n}\big)$ for $|y|\ge2^{-n}$, we obtain
\begin{align*}
2^{2n}\int_{\{|y|>2^{-n}\}}\o(|y|)\,\left|V\left(2^ny\right)\right|\,dy&\le
2\cdot2^{3n}\o\big(2^{-n}\big)\int_{\{|y|>2^{-n}\}} |y|\,\left|V\left(2^ny\right)\right|\,dy\\[.2cm]
&=2\,\o\big(2^{-n}\big)\int_{\{|y|>1\}} |y|\cdot\left|V\left(y\right)\right|\,dy
\le\const\o\big(2^{-n}\big).
\end{align*}
This proves \rf{VPn}. $\bl$

\begin{cor}
\label{Wnn}
There exists $c>0$ such that for every modulus of continuity $\o$ and for every
$f\in\L_\o\big(\R^2\big)$, the following inequalities hold:
$$
\|f*W_n\|_{L^\be}\le c\,\o\big(2^{-n}\big)\|f\|_{\L_\o(\R^2)},\quad n\in\Z.
$$
\end{cor}

\

\section{\bf Operator ideals}
\setcounter{equation}{0}
\label{Opid}

\

In this section we give a brief introduction to quasinormed ideals of operators on Hilbert space.
Recall a functional $\|\cdot\|:X\to[0,\be)$ on a vector space
$X$ is called a {\it quasinorm} on $X$ if 

(i) $\|x\|=0$ if and only if $x=\0$;

(ii) $\|\a x\|=|\a|\cdot\|x\|$, for every $x\in X$ and $\a\in\C$;

(iii) there exists a positive number $c$ such that $\|x+y\|\le c\big(\|x\|+\|y\|)$ for every $x$ and $y$ in 
$X$.

We say that a sequence $\{x_j\}_{j\ge1}$ of vectors of a {\it quasinormed space} $X$ converges to 
$x\in X$ if $\lim\limits_{j\to\be}\|x_j-x\|=0$. It is well known that there exists a translation invariant metric on $X$ which induces an equivalent topology on $X$. A quasinormed space is called {\it quasi-Banach} if it is complete.

Recall that for a bounded linear operator $T$ on Hilbert space, the singular
values $s_j(T)$, $j\ge0$, are defined by
$$
s_j(T)=\inf\big\{\|T-R\|:~\rank R\le j\big\}.
$$
Clearly, $s_0(T)=\|T\|$ and $T$ is compact if and only if $s_j(T)\to0$ as $j\to\be$. We also introduce the sequence $\{\s_n(T)\}_{n\ge0}$
defined by
\bay
\label{sin}
\s_n(T)\df\frac1{n+1}\sum_{j=0}^ns_j(T).
\ey

\medskip

{\bf Definition.} Let $\h$ be a Hilbert space and let ${\frak I}$ be a linear manifold in the set $\B(\h)$ of bounded linear operators on $\h$ that is equipped with a quasi-norm $\|\cdot\|_{\frak I}$ that makes
$\fI$ a quasi-Banach space.
We say that ${\frak I}$ is a {\it quasinormed ideal} if for every $A$ and $B$ in $\B(\h)$ and 
$T\in{\frak I}$,
\bay
\label{qni}
ATB\in{\frak I}\qm\|ATB\|_\fI\le\|A\|\cdot\|B\|\cdot\|T\|_\fI.
\ey
A quasinormed ideal $\fI$ is called a {\it normed ideal} if $\|\cdot\|_\fI$ is a norm.

Note that we do not require that $\fI\ne\B(\h)$.

\medskip

It is easy to see that if $T_1$ and $T_2$ are operators in a quasinormed ideal $\fI$ and 
$s_j(T_1)=s_j(T_2)$ for $j\ge0$, then $\|T_1\|_\fI=\|T_2\|_\fI$. Thus there exists a function 
$\Psi=\Psi_\fI$ defined on the set of nonincreasing sequences of nonnegative real numbers with values in $[0,\be]$
such that $T\in\fI$ if and only if $\Psi\big(s_0(T),s_1(T),s_2(T),\cdots~)<\be$ and
$$
\|T\|_\fI=\Psi\big(s_0(T),s_1(T),s_2(T),\cdots~),\quad T\in\fI.
$$
If $T$ is an operator from a Hilbert space $\h_1$ to a Hilbert space $\h_2$, we say that $T$ belongs to 
$\fI$ if $\Psi\big(s_0(T),s_1(T),s_2(T),\cdots~)<\be$.

For a quasinormed ideal $\fI$ and a positive number $p$, we define the quasinormed ideal $\fI^{\{p\}}$ by
$$
\fI^{\{p\}}=\left\{T:~\big(T^*T\big)^{p/2}\in\fI\right\},
\quad\|T\|_{\fI^{\{p\}}}\df\left\|(T^*T\big)^{p/2}\right\|_\fI^{1/p}.
$$

If $T$ is an operator on a Hilbert space $\h$ and $d$ is a positive integer, we denote by $[T]_d$
the operator $\bigoplus\limits_{j=1}^dT_j$ on the orthogonal sum $\bigoplus\limits_{j=1}^d\h$ of $d$ copies of $\h$, where 
$T_j=T$, $1\le j\le d$. It is easy to see that 
$$
s_n\big([T]_d\big)=s_{[n/d]}(T),\quad n\ge0,
$$
where $[x]$ denotes the largest integer that is less than or equal to $x$.

We denote by $\b_{\fI,d}$ the quasinorm of the transformer $T\mapsto[T]_d$
on $\fI$. Clearly, the sequence $\{\b_{\fI,d}\}_{d\ge1}$ is nondecreasing and {\it submultiplicative}, i.e.,
$\b_{\fI,d_1d_2}\le\b_{\fI,d_1}\b_{\fI,d_2}$. It is well known (see e.g., \S\,3 of \cite{AP3}) that the last inequality implies that
\bay
\label{lif}
\lim_{d\to\be}\frac{\log\b_{\fI,d}}{\log d}=\inf_{d\ge2}\frac{\log\b_{\fI,d}}{\log d}.
\ey

\medskip

{\bf Definition.} If $\fI$ is a quasinormed ideal, the number
$$
\b_\fI\df\lim_{d\to\be}\frac{\log\b_{\fI,d}}{\log d}=\inf_{d\ge2}\frac{\log\b_{\fI,d}}{\log d}
$$
is called the {\it upper Boyd index of} $\fI$.

\medskip

It is easy to see that $\b_\fI\le1$ for an arbitrary normed ideal $\fI$. It is also clear that $\b_\fI<1$
if and only if $\lim\limits_{d\to\be}d^{-1}\b_{\fI,d}=0$.

Note that the upper Boyd index does not change if we replace the initial quasinorm in the quasinormed ideal with an equivalent one that also satisfies \rf{qni}. It is also easy to see that
$$
\b_{\fI^{\{p\}}}=p^{-1}\b_\fI.
$$

The proof of the following fact can be found in \cite{AP3}, \S\,3.

\medskip

{\bf Theorem on ideals with upper Boyd index less than 1.}
{\em Let $\fI$ be a quasinormed ideal. The following are equivalent:

{\em(i)} $\b_\fI<1$;

{\em(ii)} for every nonincreasing sequence $\{s_n\}_{\ge0}$
of nonnegative numbers,
\bay
\label{xi}
\Psi_\fI\Big(\{\s_n\}_{n\ge0}\Big)\le\const\Psi_\fI\Big(\{s_n\}_{n\ge0}\Big),
\ey
where 
$\s_n\df(1+n)^{-1}\sum\limits_{j=0}^ns_j$}.

\medskip

For a {\it normed} ideal $\fI$ let
$\bs{C}_\fI$ be the best possible constant in inequality \rf{xi}.
Then (see \cite{AP3}, \S\,3)
\bay
\label{Ce}
\bs{C}_\fI\le3\sum_{k=0}^\be2^{-k}\b_{\fI,2^k}.
\ey

\medskip

Let $\bS_p$, $0<p<\be$, be the Schatten--von Neumann class of operators $T$ on Hilbert space
such that
$$
\|T\|_{\bS_p}\df\left(\sum_{j\ge0}\big(s_j(T)\big)^p\right)^{1/p}.
$$
This is a normed ideal for $p\ge1$. We denote by $\bS_{p,\be}$, $0<p<\be$, the ideal that consists of operators $T$ on Hilbert space such that
$$
\|T\|_{\bS_{p,\be}}\df\left(\sup_{j\ge0}(1+j)\big(s_j(T)\big)^p\right)^{1/p}.
$$
The quasinorm $\|\cdot\|_{p,\be}$ is not a norm, but it is equivalent to a norm if $p>1$.
It is easy to see that 
$$
\b_{\bS_p}=\b_{\bS_{p,\be}}=\frac1p,\quad 0<p<\be.
$$
Thus $\bS_p$ and $\bS_{p,\be}$ with $p>1$ satisfy the hypotheses of Theorem on ideals with upper Boyd index less than 1.

It follows easily from \rf{Ce} that for $p>1$,
$$
\bs{C}_{\bS_p}\le3\big(1-2^{1/p-1}\big)^{-1}.
$$

Suppose now that $\fI$ is a quasinormed ideal of operators on Hilbert space.
With a nonnegative integer $l$ we associate the ideal $^{(l)}\fI$ 
that consists of all bounded linear operators on Hilbert space and
is equipped with the norm
$$
\Psi_{^{(l)}\fI}(s_0,s_1,s_2,\cdots)=\Psi(s_0,s_1,\cdots,s_l,0,0,\cdots).
$$
It is easy to see that for every bounded operator $T$,
\begin{align*}
\|T\|_{^{(l)}\fI}&=\sup\big\{\|RT\|_\fI:~\|R\|\le1,~\rank R\le l+1\big\}
\\[.2cm]
&=\sup\big\{\|TR\|_\fI:~\|R\|\le1,~\rank R\le l+1\big\}.
\end{align*}

It is easy to verify (see \cite{AP3}, \S\,3) that if
 $\fI$ is a quasinormed ideal, then for all $l\ge0$,
\bay
\label{ClI}
\bs{C}_{^{(l)}\fI}\le\bs{C}_\fI.
\ey

Note that if $\fI=\bS_p$, $p\ge1$, then $\bS_p^l\df{^{(l)}\bS_p}$
is the normed ideal that consists of all bounded linear operators equipped with the norm
$$
\|T\|_{\bS_p^l}\df\left(\sum_{j=0}^l\big(s_j(T)\big)^p\right)^{1/p}.
$$
It is well known that $\|\cdot\|_{\bS_p^l}$ is a norm for $p\ge1$ (see \cite{BS0}).

It is also well known (see \cite{AP3}, \S\,3) that
\bay
\label{rpq}
\|T_1T_2\|_{\bS_r^l}\le\|T_1\|_{\bS_p^l}\|T_2\|_{\bS_q^l},
\ey
where $T_1$ and $T_2$ bounded operator on Hilbert space and $1/p+1/q=1/r$.

We say that a quasinormed ideal $\fI$ has {\it majorization property}  (respectively {\it weak majorization property}) if the conditions
$$
T_1\in\fI,\quad T_2\in\B,\quad
\mbox{and}\quad
\s_l(T_2)\le\s_l(T_1)\quad
\mbox{for all}\quad
l\ge0
$$
imply that
$$
T_2\in\fI\quad\mbox{and}\quad\|T_2\|_{\fI}\le\|T_1\|_{\fI}\quad(\text{respectively}\quad
\|T_2\|_{\fI}\le C\|T_1\|_{\fI})
$$
(see \cite{GK}).
Note that if a quasinormed ideal $\fI$ has weak majorization property, then we can introduce on it the following new equivalent quasinorm:
$$
\|T\|_{\widetilde\fI}\df\sup\{\|R\|_{\fI}:~\s_l(R)\le\s_l(T)\,\,\,\text{for all}\,\,\,l\ge0\}
$$
such that
$(\fI,\|\cdot\|_{\widetilde\fI})$ has majorization property.

It is well known that every separable normed ideal and every normed ideal
that is dual to a separable normed ideal has majorization property, see \cite{GK}.
Clearly, $\bS_1\subset\fI$ for every quasinormed ideal $\fI$
with majorization property.
Note also that every quasinormed ideal $\fI$ with $\b_{\fI}<1$ has weak majorization property (see, for example, \S\,3 of \cite{AP3} and
\S\,3 of \cite{AP4}).

We need the following fact on interpolation properties of quasinormed ideals that have majorization property (see e.g., \cite{AP4}):

\medskip

{\bf Theorem on interpolation of quasinormed ideals.}
{\it Let $\fI$ be a quasinormed ideal with majorization property
and let $\frak A:\frak L\to\frak L$ be a linear transformer
on a linear subset $\frak L$ of $\B$ such that $\frak L\cap\bS_1$
is dense in $\bS_1$. Suppose that
$\|\frak A T\|\le\|T\|$ and $\|\frak A T\|_{\bS_1}\le\|T\|_{\bS_1}$
for all $T\in\frak L$. Then $\|\frak A T\|_{\fI}\le\|T\|_{\fI}$
for every  $T\in\frak L$}.

We refer the reader to \cite{GK} and \cite{BS0} for further information on singular values and normed ideals of operators on Hilbert space.

\

\section{\bf Double operator integrals}
\setcounter{equation}{0}
\label{dois}

\

In this subsection we give a brief introduction in double  operator integrals. Double operator integrals appeared in the paper \cite{DK} by Daletskii and S.G. Krein. However, the beautiful theory of double operator integrals was developed later by Birman and Solomyak in \cite{BS1}, \cite{BS2}, and \cite{BS3}, see also their survey \cite{BS4}.

Let $(\X,E_1)$ and $(\Y,E_2)$ be spaces with spectral measures $E_1$ and $E_2$
on a Hilbert space $\h$. The idea of Birman and Solomyak is to define first
double operator integrals
\bay
\label{doi}
\int\limits_\X\int\limits_\Y\Phi(x,y)\,d E_1(x)T\,dE_2(y),
\ey
for bounded measurable functions $\Phi$ and operators $T$
of Hilbert Schmidt class $\bS_2$. Consider the spectral measure $\E$ whose values are orthogonal
projections on the Hilbert space $\bS_2$, which is defined by
$$
\E(\L\times\D)T=E_1(\L)TE_2(\D),\quad T\in\bS_2,
$$ 
$\L$ and $\D$ being measurable subsets of $\X$ and $\Y$. It was shown in \cite{BS} that $\E$ extends to a spectral measure on 
$\X\times\Y$ and if $\Phi$ is a bounded measurable function on $\X\times\Y$, by definition,
$$
\int\limits_\X\int\limits_\Y\Phi(x,y)\,d E_1(x)T\,dE_2(y)=
\left(\,\,\int\limits_{\X\times\Y}\Phi\,d\E\right)T.
$$
Clearly,
$$
\left\|\int\limits_\X\int\limits_\Y\Phi(x,y)\,dE_1(x)T\,dE_2(y)\right\|_{\bS_2}
\le\|\Phi\|_{L^\be}\|T\|_{\bS_2}.
$$
If 
$$
\int\limits_\X\int\limits_\Y\Phi(x,y)\,d E_1(x)T\,dE_2(y)\in\bS_1
$$
for every $T\in\bS_1$, we say that $\Phi$ is a {\it Schur multiplier of $\bS_1$ associated with 
the spectral measures $E_1$ and $E_2$}. 

In this case the transformer
\bay
\label{tra}
T\mapsto\int\limits_{\Y}\int\limits_{\X}\Phi(x,y)\,d E_2(y)\,T\,dE_1(x),\quad T\in \bS_2,
\ey
extends by duality to a bounded linear transformer on the space of bounded linear operators on $\h$
and we say that the function $\Psi$ on $\Y\times\X$ defined by 
$$
\Psi(y,x)=\Phi(x,y)
$$
is {\it a Schur multiplier (with respect to $E_2$ and $E_1$) of the space of bounded linear operators}.
We denote the space of such Schur multipliers by $\fM(E_2,E_1)$.
The norm of $\Psi$ in $\fM(E_2,E_1)$ is, by definition, the norm of the
transformer \rf{tra} on the space of bounded linear operators.

In \cite{BS3} it was shown that if $A$ and $B$ are a self-adjoint operators (not necessarily bounded) such that $A-B$ is bounded
 and if $f$ is a continuously differentiable 
function on $\R$ such that the divided difference $\dg f$,
$$
\big(\dg f\big)(x,y)=\frac{f(x)-f(y)}{x-y},
$$
is a Schur multiplier
of $\bS_1$ with respect to the spectral measures of $A$ and $B$, then
\bay
\label{BSF}
f(A)-f(B)=\iint\big(\dg f\big)(x,y)\,dE_{A}(x)(A-B)\,dE_B(y)
\ey
and
$$
\|f(A)-f(B)\|\le\const\|f\|_{\fM(E_A,E_{B})}\|A-B\|,
$$
i.e., {\it $f$ is an operator Lipschitz function}.

It is easy to see that if a function $\Phi$ on $\X\times\Y$ belongs to the {\it projective tensor
product}
$L^\be(E_1)\hat\otimes L^\be(E_2)$ of $L^\be(E_1)$ and $L^\be(E_2)$ (i.e., $\Phi$ admits a representation
$$
\Phi(x,y)=\sum_{n\ge0}\f_n(x)\psi_n(y),
$$
where $\f_n\in L^\be(E_1)$, $\psi_n\in L^\be(E_2)$, and
$$
\sum_{n\ge0}\|\f_n\|_{L^\be}\|\psi_n\|_{L^\be}<\be),
$$
then $\Phi\in\fM(E_1,E_2)$.
For such functions $\Phi$ we have
$$
\int\limits_\X\int\limits_\Y\Phi(x,y)\,dE_1(x)T\,dE_2(y)=
\sum_{n\ge0}\left(\,\int\limits_\X\f_n\,dE_1\right)T\left(\,\int\limits_\Y\psi_n\,dE_2\right).
$$ 

More generally, $\Phi\in\fM(E_1,E_2)$ if $\Phi$ 
belongs to the {\it integral projective tensor product} $L^\be(E_1)\hat\otimes_{\rm i}
L^\be(E_2)$ of $L^\be(E_1)$ and $L^\be(E_2)$, i.e., $\Phi$ admits a representation
\bay
\label{ipt}
\Phi(x,y)=\int_\O \f(x,w)\psi(y,w)\,d\l(w),
\ey
where $(\O,\l)$ is a $\s$-finite measure space, $\f$ is a measurable function on $\X\times \O$,
$\psi$ is a measurable function on $\Y\times \O$, and
\bay
\label{ir}
\int_\O\|\f(\cdot,w)\|_{L^\be(E_1)}\|\psi(\cdot,w)\|_{L^\be(E_2)}\,d\l(w)<\be.
\ey
If $\Phi\in L^\be(E_1)\hat\otimes_{\rm i}L^\be(E_2)$, then
$$
\int\limits_\X\int\limits_\Y\Phi(x,y)\,dE_1(x)T\,dE_2(y)=
\int\limits_\O\left(\,\int\limits_\X\f(x,w)\,dE_1(x)\right)T
\left(\,\int\limits_\Y\psi(y,w)\,dE_2(y)\right)\,d\l(w).
$$
Clearly, the function 
$$
s\mapsto \left(\,\int_\X\f(x,w)\,dE_1(x)\right)T
\left(\,\int_\Y\psi(y,w)\,dE_2(y)\right)
$$ 
is weakly measurable and
$$
\int\limits_\O\left\|\left(\,\int\limits_\X\f(x,s)\,dE_1(x)\right)T
\left(\,\int\limits_\Y\psi(y,w)\,dE_2(w)\right)\right\|\,d\l(w)<\be.
$$

It turns out that all Schur multipliers can be obtained in this way. More precisely, the following
result holds (see \cite{Pe2}):

\medskip

{\bf Theorem on Schur multipliers.} {\em Let $\Phi$ be a measurable function on 
$\X\times\Y$. The following are equivalent:

{\rm (i)} $\Phi\in\fM(E_1,E_2)$;

{\rm (ii)} $\Phi\in L^\be(E_1)\hat\otimes_{\rm i}L^\be(E_2)$;

{\rm (iii)} there exist measurable functions $\f$ on $\X\times\O$ and $\psi$ on $\Y\times\O$ such that
{\em\rf{ipt}} holds and
\bay
\label{bs}
\left\|\left(\int_\O|\f(\cdot,w)|^2\,d\l(w)\right)^{1/2}\right\|_{L^\be(E)}
\left\|\left(\int_\O|\psi(\cdot,w)|^2\,d\l(w)\right)^{1/2}\right\|_{L^\be(F)}<\be.
\ey
}

The implication (iii)$\imp$(i) was established in \cite{BS3}. 
In the case of matrix Schur multipliers (this corresponds to discrete spectral measures
of multiplicity 1) the fact that (i) implies (ii) was proved in \cite{Be}.

Note that the infimum of the left-hand side in \rf{bs} over all representations of the form \rf{ipt} is the so-called Haagerup tensor norm of two $L^\be$ spaces.

It is interesting to observe that if $\f$ and $\psi$ satisfy \rf{ir}, then they also satisfy
\rf{bs}, but the converse is false. However, if $\Phi$ admits a representation of the form \rf{ipt}
with $\f$ and $\psi$ satisfying \rf{bs}, then it also admits a (possibly different) representation of the
form \rf{ipt} with $\f$ and $\psi$ satisfying \rf{ir}.
We refer the reader to \cite{Pi} for related problems.

It is also well known that $\fM(E_1,E_2)$ is a Banach algebra (see \cite{Pe2}).

To conclude this section, we would like to observe that it follows from  
the Theorem on interpolation of quasinormed ideals (see \S\,\ref{Opid}) that if $\Phi\in\fM(E_1,E_2)$ and $\fI$ is a quasinormed ideal with majorization property, then
$$
T\in\fI\quad\Longrightarrow\quad
\int\limits_{\X}\int\limits_{\Y}\Phi(x,y)\,d E_1(x)\,T\,dE_2(y)\in\fI
$$
and
\bay
\label{intn}
\left\|\int\limits_{\X}\int\limits_{\Y}\Phi(x,y)\,d E_1(x)\,T\,dE_2(y)\right\|
\le\|\Phi\|_{\fM(E_1,E_2)}\|T\|_\fI.
\ey

\

\section{\bf The basic formula in terms of double operator integrals}
\setcounter{equation}{0}
\label{bfdoi}

\

Recall that a function $f$ on $\R^2$ is called {\it operator Lipschitz} if 
\bay
\label{OL2}
\|f(N_1)-f(N_2)\|\le\const\|N_1-N_2\|
\ey
for every normal operators $N_1$ and $N_2$ on Hilbert space.
Clearly, if $f$ is operator Lipschitz, then $f$ is a Lipschitz function. The converse is false, because it is false for self-adjoint operators (see the Introduction).

%We show in \S\,+++ that if $f$ is operator Lipschitz, then inequality \rf{OL2} holds for arbitrary (not necessarily bounded) normal operators $N_1$ and $N_2$ whose difference $N_1-N_2$ is bounded (the latter means that $N_1$ and $N_2$ have the same domains and the difference is bounded on the common domain).

The first natural try to prove that a function on $\R^2$ is operator Lipschitz is to attempt to generalize formula \rf{BSF} to the case of normal operators.
Suppose that the divided difference
$$
(z_1,z_2)\mapsto\frac{f(z_1)-f(z_2)}{z_1-z_2},\quad z_1,~z_2\in\C,
$$
is a Schur multiplier with respect to arbitrary Borel spectral measures on $\C$. Then as in the case of self-adjoint operators, 
for arbitrary normal operators $N_1$ and $N_2$ with bounded difference $N_1-N_2$, the following formula holds
\bay
\label{ddn}
f(N_1)-f(N_2)=\iint\limits_{\C\times\C}\frac{f(z_1)-f(z_2)}{z_1-z_2}
\,dE_1(z_1)(N_1-N_2)\,dE_2(z_2),
\ey 
where $E_j$ is the spectral measure of $N_i$, $i=1,\,2$. Moreover, in this case $f$ is operator Lipschitz.

However, it follows from the results of \cite{JW} that under the above assumptions $f$ must have complex derivative everywhere. In other words, $f$ must be an entire function. In addition to this $f$ must be Lipschitz. Therefore in this case $f$ is a linear function, but the fact that linear functions are operator Lipschitz is obvious.

Thus to prove that a given function on $\R^2$ is operator Lipschitz, we have to find something different.

To state the main results of this section, we introduce the following notation. Given normal operators $N_1$ and $N_2$ on Hilbert space, we put
$$
A_j\df\re N_j,\quad B_j\df\im N_j,\quad E_j\quad\mbox{is the spectral measure of}\quad N_j,\quad j=1,\,2.
$$
In other words, $N_j=A_j+{\rm i}B_j$, $j=1,\,2$, where $A_j$ and $B_j$ are self-adjoint operators. Since the operators $N_j$ are normal, $A_j$ commutes with $B_j$. 

With a function $f$ on $\R^2$ that has partial derivatives everywhere, we associate the following divided differences
$$
\big(\dg_xf\big)(z_1,z_2)\df\frac{f(x_1,y_2)-f(x_2,y_2)}{x_1-x_2},
\quad z_1,\,z_2\in\C.
$$
and
$$
\big(\dg_yf\big)(z_1,z_2)\df\frac{f(x_1,y_1)-f(x_1,y_2)}{y_1-y_2},
\quad z_1,\,z_2\in\C.
$$
Throughout the paper we use the notation
$$
x_j\df\re z_j,\quad y_j\df\im z_j,\quad  j=1,\,2.
$$
Note that in the above definition by the values of $\dg_xf$ and $\dg_yf$ 
on the sets 
$$
\{(z_1,z_2):~x_1=x_2\}\quad\mbox{and}\quad\{(z_1,z_2):~y_1=y_2\}
$$
we mean the corresponding partial derivatives of $f$.

Let us now state the main results of this section.

\begin{thm}
\label{SMdd}
Let $f$ be a continuous bounded function on $\R^2$ whose Fourier transform $\F f$ has compact support. Then the functions $\dg_xf$ and $\dg_yf$ are Schur multipliers with respect to arbitrary Borel spectral measures $E_1$ and $E_2$. 

Moreover, if
$$
\supp\F f\subset\{\z\in\C:~|\z|\le\s\},\quad\s>0,
$$ 
then
\bay
\label{nerxy}
\|\dg_xf\|_{\fM(E_1,E_2)}\le\const\s\|f\|_{L^\be}
\quad\mbox{and}\quad\|\dg_yf\|_{\fM(E_1,E_2)}\le \const\s\|f\|_{L^\be}.
\ey
\end{thm}

\begin{thm}
\label{doin}
Let $f$ be a continuous bounded function on $\R^2$ whose Fourier transform $\F f$ has compact support. Suppose that $N_1$ and $N_2$ are normal operators such that the operator $N_1-N_2$ is bounded. Then
\begin{align}
\label{if}
f(N_1)-f(N_2)&=\iint\limits_{\C^2}\big(\dg_yf\big)(z_1,z_2)\,
dE_1(z_1)(B_1-B_2)\,dE_2(z_2)\nonumber\\[.2cm]
&+\iint\limits_{\C^2}\big(\dg_xf\big)(z_1,z_2)\,
dE_1(z_1)(A_1-A_2)\,dE_2(z_2).
\end{align}
\end{thm}

We postpone the proof of Theorem \ref{SMdd} till the next section. Let us deduce here Theorem \ref{doin} from Theorem \ref{SMdd}.  

\medskip

{\bf Proof of Theorem \ref{doin}.} Consider first the case when $N_1$ and $N_2$ are bounded operators. Put 
$$
d=\max\big\{\|N_1\|,\|N_2\|\big\}\quad\mbox{and}\quad D\df\{\z\in\C:~|\z|\le d\}.
$$

By Theorem \ref{SMdd}, both $\dg_yf$ and $\dg_xf$ are Schur multipliers.
We have
\begin{align*}
\iint\limits_{\C^2}&\big(\dg_yf\big)(z_1,z_2)\,
dE_1(z_1)(B_1-B_2)\,dE_2(z_2)\\[.2cm]
&=\iint\limits_{D\times D}\big(\dg_yf\big)(z_1,z_2)\,
dE_1(z_1)(B_1-B_2)\,dE_2(z_2)\\[.2cm]
&=\iint\limits_{D\times D}\big(\dg_yf\big)(z_1,z_2)\,dE_1(z_1)B_1\,dE_2(z_2)
-\iint\limits_{D\times D}\big(\dg_yf\big)(z_1,z_2)\,dE_1(z_1)B_2\,dE_2(z_2)
\\[.2cm]
&=\iint\limits_{D\times D}y_1\big(\dg_yf\big)(z_1,z_2)\,dE_1(z_1)\,dE_2(z_2)
-\iint\limits_{D\times D}y_2\big(\dg_yf\big)(z_1,z_2)\,dE_1(z_1)\,dE_2(z_2)
\\[.2cm]
&=\iint\limits_{D\times D}(y_1-y_2)\big(\dg_yf\big)(z_1,z_2)\,dE_1(z_1)\,dE_2(z_2)\\[.2cm]
&=\iint\limits_{D\times D}\big(f(x_1,y_1)-f(x_1,y_2)\big)\,dE_1(z_1)\,dE_2(z_2).
\end{align*}
Since $\fM(E_1,E_2)$ is a Banach algebra, it is easy to see that the function
$$
(z_1,z_2)\mapsto f(x_1,y_1)-f(x_1,y_2)=(y_1-y_2)\big(\dg_yf\big)(z_1,z_2)
$$
is a Schur multiplier.
Similarly,
$$
\iint\limits_{\C^2}\!\big(\dg_xf\big)(z_1,z_2)\,
dE_1(z_1)(A_1-A_2)\,dE_2(z_2)
=\iint\limits_{D\times D}\!\big(f(x_1,y_2)-f(x_2,y_2)\big)\,dE_1(z_1)\,dE_2(z_2).
$$
It follows that
\begin{align*}
\iint\limits_{\C^2}&\big(\dg_yf\big)(z_1,z_2)\,
dE_1(z_1)(B_1-B_2)\,dE_2(z_2)\\[.2cm]
+&\iint\limits_{\C^2}\big(\dg_xf\big)(z_1,z_2)\,
dE_1(z_1)(A_1-A_2)\,dE_2(z_2)\\[.2cm]
&=\iint\limits_{D\times D}\big(f(x_1,y_1)-f(x_2,y_2)\big)\,dE_1(z_1)\,dE_2(z_2)\\[.2cm]
&=\iint\limits_{D\times D}f(x_1,y_1)\,dE_1(z_1)\,dE_2(z_2)
-\iint\limits_{D\times D}f(x_2,y_2)\,dE_1(z_1)\,dE_2(z_2)\\[.2cm]
&=f(N_1)-f(N_2).
\end{align*} 

Consider now the case when $N_1$ and $N_2$ are unbounded. Put
$$
P_k\df E_1\big(\{\z\in\C:|\z|\le k\}\big)\quad\mbox{and}\quad
Q_k\df E_2\big(\{\z\in\C:|\z|\le k\}\big),\quad k>0.
$$
Then
$$
N_{1,k}\df P_kN_1\quad\mbox{and}\quad N_{2,k}\df Q_kN_2
$$
are bounded normal operators. Denote by $E_{j,k}$ the spectral measure of
$N_{j,k}$, $j=1,\,2$. It is easy to see that
$$
N_{1,k}=P_kA_1+{\rm i}P_kB_1,\quad \mbox{and}\quad
N_{2,k}=A_2Q_k+{\rm i}B_2Q_k,\quad k>0.
$$

%Clearly,
%$$
%\lim_{k\to\be}P_k\big(f(N_1)-f(N_2)\big)Q_k=f(N_1)-f(N_2),
%$$
%\begin{align*}
%\lim_{k\to\be}P_k&\left(\iint\limits_{\C^2}\big(\dg_yf\big)(z_1,z_2)\,
%dE_1(z_1)(B_1-B_2)\,dE_2(z_2)\right)Q_k\\[.2cm]
%&=
%\iint\limits_{\C^2}\big(\dg_yf\big)(z_1,z_2)\,
%dE_1(z_1)(B_1-B_2)\,dE_2(z_2),
%\end{align*}
%and
%\begin{align*}
%\lim_{k\to\be}P_k&\left(\iint\limits_{\C^2}\big(\dg_xf\big)(z_1,z_2)\,
%dE_1(z_1)(A_1-A_2)\,dE_2(z_2)\right)Q_k\\[.2cm]
%&=
%\iint\limits_{\C^2}\big(\dg_xf\big)(z_1,z_2)\,
%dE_1(z_1)(A_1-A_2)\,dE_2(z_2).
%\end{align*}
%in the strong operator topology.

We have
\begin{align*}
P_k&\left(\iint\limits_{\C^2}\big(\dg_yf\big)(z_1,z_2)\,
dE_1(z_1)(B_1-B_2)\,dE_2(z_2)\right)Q_k\\[.2cm]
&=P_k\left(\iint\limits_{\C^2}\big(\dg_yf\big)(z_1,z_2)\,
dE_{1,k}(z_1)(P_kB_1-B_2Q_k)\,dE_{2,k}(z_2)\right)Q_k
\end{align*}
and
\begin{align*}
P_k&\left(\iint\limits_{\C^2}\big(\dg_xf\big)(z_1,z_2)\,
dE_1(z_1)(A_1-A_2)\,dE_2(z_2)\right)Q_k\\[.2cm]
&=P_k\left(\iint\limits_{\C^2}\big(\dg_yf\big)(z_1,z_2)\,
dE_{1,k}(z_1)(P_kA_1-A_2Q_k)\,dE_{2,k}(z_2)\right)Q_k.
\end{align*}

If we apply identity \rf{if} to the bounded normal operators $N_{1,k}$ and $N_{2,k}$, we obtain
\begin{align*}
P_k\big(f(N_{1,k})&-f(N_{2,k})\big)Q_k=\\[.2cm]
=&P_k\left(\iint\limits_{\C^2}\big(\dg_yf\big)(z_1,z_2)\,
dE_{1,k}(z_1)(P_kB_1-B_2Q_k)\,dE_{2,k}(z_2)\right)Q_k\\[.2cm]
&+P_k\left(\iint\limits_{\C^2}\big(\dg_yf\big)(z_1,z_2)\,
dE_{1,k}(z_1)(P_kA_1-A_2Q_k)\,dE_{2,k}(z_2)\right)Q_k.
\end{align*}
Since obviously,
$$
P_k\big(f(N_{1,k})-f(N_{2,k})\big)Q_k=P_k\big(f(N_{1})-f(N_{2})\big)Q_k,
$$
we have
\begin{align*}
P_k\big(f(N_{1})&-f(N_{2})\big)Q_k=\\[.2cm]
=&P_k\left(\iint\limits_{\C^2}\big(\dg_yf\big)(z_1,z_2)\,
dE_1(z_1)(B_1-B_2)\,dE_2(z_2)\right)Q_k\\[.2cm]
&+P_k\left(\iint\limits_{\C^2}\big(\dg_xf\big)(z_1,z_2)\,
dE_1(z_1)(A_1-A_2)\,dE_2(z_2)\right)Q_k.
\end{align*}
It remains to pass to the limit in the strong operator topology. $\bl$

We would like to extend formula \rf{if} to the case of arbitrary functions $f$ in $B_{\be1}^1\big(\R^2\big)$. Since $B_{\be1}^1\big(\R^2\big)$ consists of Lipschitz functions, it follows that for $f\in B_{\be1}^1(\R^2)$, 
\bay
\label{1+z}
|f(\z)|\le\const(1+|\z|),\quad\z\in\C.
\ey

Hence, for $f\in B_{\be1}^1(\R^2)$,
$$
D_{f(N)}\supset D_N.
$$

\begin{thm}
\label{dlya}
Let $N_1$ and $N_2$ be normal operators such that $N_1-N_2$ is bounded.
Then {\em\rf{if}} holds for every $f\in B_{\be1}^1\big(\R^2\big)$.
\end{thm}

\Pf It suffices to prove that for $u\in D_{N_1}=D_{N_2}$, 
\begin{align*}
\big(f(N_1)-f(N_2)\big)u&=\left(\iint\limits_{\C^2}\big(\dg_yf\big)(z_1,z_2)\,
dE_1(z_1)(B_1-B_2)\,dE_2(z_2)\right)u\\[.2cm]
&+\left(\iint\limits_{\C^2}\big(\dg_xf\big)(z_1,z_2)\,
dE_1(z_1)(A_1-A_2)\,dE_2(z_2)\right)u.
\end{align*}
Indeed, if $N$ is a normal operator and $f$ satisfies \rf{1+z},
then $f(N)$ is the closure of its restriction to the domain of $N$.

We have
$$
\big(f(N_1)-f(N_2)\big)u=\big(\big(f-f(0)\big)(N_1)\big)u-
\big(\big(f-f(0)\big)(N_2)\big)u,
$$
\bay
\label{N1u}
\big(\big(f-f(0)\big)(N_1)\big)u=
\sum_{n\in\Z}\big(\big(f_n-f_n(0)\big)(N_1)\big)u,
\ey
and
\bay
\label{N2u}
\big(\big(f-f(0)\big)(N_2)\big)u=
\sum_{n\in\Z}\big(\big(f_n-f_n(0)\big)(N_2)\big)u,
\ey
where the functions $f_n$ are defined by \rf{fn}.
Moreover, the series on the right-hand sides of \rf{N1u} and \rf{N2u} converge absolutely in the norm.

Thus 
$$
\big(f(N_1)-f(N_2)\big)u=
\sum_{n\in\Z}\big(f_n(N_1)-f_n(N_2)\big)u.
$$
It remains to observe that
\begin{align*}
&\iint\limits_{\C^2}\big(\dg_yf\big)(z_1,z_2)\,dE_1(z_1)(B_1-B_2)\,dE_2(z_2)
\\[.2cm]
=&\sum_{n\in\Z}
\iint\limits_{\C^2}\big(\dg_yf_n\big)(z_1,z_2)\,dE_1(z_1)(B_1-B_2)\,dE_2(z_2)
\end{align*}
and
\begin{align*}
&\iint\limits_{\C^2}\big(\dg_xf\big)(z_1,z_2)\,dE_1(z_1)(A_1-A_2)\,dE_2(z_2)
\\[.2cm]
=&\sum_{n\in\Z}
\iint\limits_{\C^2}\big(\dg_xf_n\big)(z_1,z_2)\,dE_1(z_1)(A_1-A_2)\,dE_2(z_2),
\end{align*}
and the series on the right-hand sides converge absolutely in the norm
which is an immediate consequence of inequalities \rf{nerxy}.  $\bl$

\

\section{\bf Proof of Theorem \ref{SMdd}}
\setcounter{equation}{0}
\label{PoT}

\

In this section we are going to prove Theorem \ref{SMdd} that gives sharp estimates for the norms of $\dg_xf$ and $\dg_yf$ in the space of Schur multipliers. Consider the function $\dg_xf$,
$$
\big(\dg_xf\big)(z_1,z_2)=\frac{f(x_1,y_2)-f(x_2,y_2)}{x_1-x_2},\quad z_1,z_2\in\C.
$$
The first natural thought would be to fix the variable $y_2$ and represent the function
$$
(x_1,x_2)\mapsto\frac{f(x_1,y_2)-f(x_2,y_2)}{x_1-x_2}
$$
in terms of the integral projective tensor product 
$L^\be\hat\otimes_{\rm i}L^\be$ in the same was as it was done in \cite{Pe3} for functions of one variable. However, it turns out that if we do this, we obtain in the integral tensor representation terms that depend on the mixed variables $(x_1,y_2)$, and so this would not help us.

The first proof of Theorem \ref{SMdd} we have found was based on a modification of the integral tensor representation obtained in \cite{Pe3} and an estimate in terms of the tensor norm \rf{bs} rather than the integral projective tensor norm. 

In this section we give a different approach based on an expansion of entire functions of exponential type $\s$ in the series in the orthogonal basis 
$\left\{\dfrac{\sin \s x}{\s x-\pi n}\right\}_{n\in\Z}$.

For a topological space $\X$, we denote by $C_{\rm b}(\X)$ the set of
bounded continuous (complex) functions on $\X$.
If $\X$ and $\Y$ are topological spaces, we denote by
$C_{\rm b}(\X)\hat\otimes_{\rm h}C_{\rm b}(\Y)$
the set of functions $\Phi$ 
on $\X\times\Y$ that admit a representation 
\bay
\label{phipsi}
\Phi(x,y)=\sum_{n\ge0}\f_n(x)\psi_n(y),\quad (x,y)\in\X\times\Y
\ey
such that $\f_n\in C_{\rm b}(\X)$, $\psi_n\in C_{\rm b}(Y)$ and
\bay
\label{phipsic}
\left(\sup_{x\in\X}\sum_{n\ge0}|\f_n(x)|^2\right)^{1/2}
%\Big(\sup_{y\in\Y}\sum_{n\ge0}|\psi_n(y)|^2\Big)^{1/2}<\be
\left(\sup_{y\in\Y}\sum_{n\ge0}|\psi_n(y)|^2\right)^{1/2}<\be.
\ey
For $\Phi\in C_{\rm b}(\X)\hat\otimes_{\rm h}C_{\rm b}(\Y)$,  its norm in $C_{\rm b}(\X)\hat\otimes_{\rm h}C_{\rm b}(\Y)$
is, by definition, the infimum of the  left-hand side of \rf{phipsic} over all
representations \rf{phipsi}.

For $\s>0$, we denote by $\mathscr E_\s$ the set of entire functions (of one complex variable)
of exponential type at most $\s$.

It follows from the results of \cite{Pe3} that
\bay
\label{m1}
f\in\mathscr E_\s\cap L^\be(\R)\quad\Longrightarrow\quad
\left\|\frac{f(x)-f(y)}{x-y}\right\|_{\frak M(E_1,E_2)}\le\const\s\|f\|_{L^\be(\R)}
\ey
for every Borel spectral measures $E_1$ and $E_2$ on $\R$.

It was shown in \cite{AP4} that inequality \rf{m1} holds with constant equal to 1.

The following result allows us to obtain an explicit representation of the divided difference $\dfrac{f(x)-f(y)}{x-y}$ as an element of 
$C_{\rm b}(\R)\hat\otimes_{\rm h}C_{\rm b}(\R)$.

\begin{thm}
\label{to}
Let $f\in\mathscr E_\s\cap L^\be(\R)$. Then
\begin{align}
\label{haa1}
\frac{f(x)-f(y)}{x-y}&=\sum_{n\in\Z}(-1)^{n}\s\cdot\frac{f(x)-f(\pi n\s^{-1})}{\s x-\pi n}
\cdot\frac{\sin \s y}{\s y-\pi n}\\[.2cm]
\label{haa2}
&=\frac1\pi\int_\R\frac{f(x)-f(t)}{x-t}\cdot\frac{\sin(\s (y-t))}{y-t}\,dt,
\quad x,~y\in\R.
\end{align}
Moreover,
\bay
\label{2vy}
\sum_{n\in\Z}\frac{|f(x)-f(\pi n\s^{-1})|^2}{(\s x-\pi n)^2}
=\frac1{\pi\s}\int_\R\frac{|f(x)-f(t)|^2}{(x-t)^2}dt\le 3\|f\|_{L^\be(\R)}^2,
\quad x\in\R,
\ey
and
\bay
\label{hiz}
\sum_{n\in\Z}\frac{\sin^2 \s y}{(\s y-\pi n)^2}=1=\frac1{\pi\s}\int_\R\frac{\sin^2(\s (y-t))}{(y-t)^2}\,dt,\quad y\in\R.
\ey
\end{thm}

\Pf Clearly, it suffices to consider the case $\s=1$. Let us first observe that the identities in \rf{hiz} are elementary and well known.

We are going to use the well-known fact that the family
$\left\{\dfrac{\sin z}{z-\pi n}\right\}_{n\in\Z}$
forms an orthogonal basis in the space $\mathscr E_1\cap L^2(\R)$,
\bay
\label{Fz}
F(z)=\sum_{n\in\Z}(-1)^{n}F(\pi n)\frac{\sin z}{z-\pi n},
\ey
and
\bay
\label{Fz1}
\sum_{n\in\Z}|F(\pi n)|^2=\frac1\pi\int_\R|F(t)|^2\,dt.
\ey
for every $F\in\mathscr E_1\cap L^2(\R)$, see, e.g., \cite{L}, Lect. 20.2, Th. 1. It follows immediately from \rf{Fz1} that
\bay
\label{Fz2}
\sum_{n\in\Z}F(\pi n)\ov{G(\pi n)}=\frac1\pi\int_\R F(t)\ov{G(t)}\,dt
\quad\mbox{for every}\quad F,~G\in\mathscr E_1\cap L^2(\R).
\ey

Given $x\in\R$, we consider the function $F$ defined by $F(\l)=\dfrac{f(x)-f(\l)}{x-\l}$, $\l\in\C$. Clearly, $F\in\mathscr E_1\cap L^2(\R)$.

It is easy to see that \rf{haa1} is a consequence of 
\rf{Fz} and the equality in \rf{2vy} is a consequence of \rf{Fz1}.
It is also easy to see that \rf{haa2} follows from \rf{Fz2}.

 It remains to prove that
$$
\frac1{\pi}\int_\R\frac{|f(x)-f(t)|^2}{(x-t)^2}dt\le 3\|f\|_{L^\be(\R)}^2
$$
for every $f\in\mathscr E_1\cap L^\be(\R)$ and $x\in\R$. Without loss of generality we may assume that
$\|f\|_{L^\be(\R)}=1$. Then $\|f^\prime\|_{L^\be(\R)}\le1$ by the Bernstein
inequality. Hence, $|f(x)-f(t)|\le\min(2,|x-t|)$, and we have
\bey
\frac1{\pi}\int_\R\frac{|f(x)-f(t)|^2}{(x-t)^2}dt\le\frac1{\pi}\int_\R\frac{\min(4,(x-t)^2)}{(x-t)^2}dt
=\frac2\pi\int_0^2dt+\frac8\pi\int_2^\be\frac{dt}{t^2}=\frac8\pi<3.\quad\bl
\eey

\medskip

{\bf Remark.} Note that the equality 
$$
\frac{f(x)-f(y)}{x-y}=\frac1\pi\int_\R\frac{f(x)-f(t)}{x-t}\cdot\frac{\sin(\s (y-t))}{y-t}\,dt
$$
is an immediate consequence of the well-known fact that
$\dfrac{\sin(\s(x-y))}{\pi(x-y)}$ is the reproducing kernel for the
functional Hilbert space $\mathscr E_1\cap L^2(\R)$.

\begin{thm}
\label{fs}
Let $\s>0$ and let $f$ be a function in $C_{\rm b}(\R^2)$ such that
$$
\supp\F f\subset\{\z\in\C:~|\z|\le\s\}.
$$
Then $\dg_xf,\,\dg_yf\in C_{\rm b}(\C)\hat\otimes_{\rm h}C_{\rm b}(\C)$,
$$
\|\dg_xf\|_{C_{\rm b}(\C)\hat\otimes_{\rm h}C_{\rm b}(\C)}\le\sqrt 3\,\s\|f\|_{L^\be(\C)}
$$
and
$$
\|\dg_yf\|_{C_{\rm b}(\C)\hat\otimes_{\rm h}C_{\rm b}(\C)}\le\sqrt 3\,\s\|f\|_{L^\be(\C)}.
$$
\end{thm}

\Pf Clearly, $f$ is the restriction to $\R^2$ of an entire function of two complex variables.
Moreover, $f(\cdot,a),\,f(a,\cdot)\in \mathscr E_\s\cap L^\be(\R)$ for every $a\in\R$.
It suffices to consider the case $\s=1$. By Theorem \ref{to}, we have
\bey
\big(\dg_xf\big)(z_1,z_2)\df\frac{f(x_1,y_2)-f(x_2,y_2)}{x_1-x_2}=
\sum_{n\in\Z}(-1)^n\frac{f(\pi n,y_2)-f(x_2,y_2)}{\pi n-x_2}\cdot\frac{\sin x_1}{x_1-\pi n}
\eey
and

\bey
\big(\dg_yf\big)(z_1,z_2)\df\frac{f(x_1,y_1)-f(x_1,y_2)}{y_1-y_2}=\sum_{n\in\Z}(-1)^n\frac{f(x_1,y_1)-f(x_1,\pi n)}{y_1-\pi n}\cdot\frac{\sin y_2}{y_2-\pi n}.
\eey

Note that the functions $\dfrac{\sin x_1}{x_1-\pi n}$ and $\dfrac{f(x_1,y_1)-f(x_1,\pi n)}{y_1-\pi n}$
depend on $z_1=(x_1,y_1)$ and do not depend on $z_2=(x_2,y_2)$ while
the functions $\dfrac{f(\pi n,y_2)-f(x_2,y_2)}{\pi n-x_2}$ and $\dfrac{\sin y_2}{y_2-\pi n}$
depend on $z_2=(x_2,y_2)$ and do not depend on $z_1=(x_1,y_1)$. Moreover, by Theorem \ref{to}
we have
\bey
\sum_{n\in\Z}\frac{|f(x_1,y_1)-f(x_1,\pi n)|^2}{(y_1-\pi n)^2}\le3\|f(x_1,\cdot)\|_{L^\be(\R)}^2
\le3\|f\|_{L^\be(\C)}^2,\\
\sum_{n\in\Z}\frac{|f(\pi n,y_2)-f(x_2,y_2)|^2}{(\pi n-x_2)^2}\le3\|f(\cdot,y_2)\|_{L^\be(\R)}^2
\le3\|f\|_{L^\be(\C)}^2,\\
\eey
and
$$
\sum_{n\in\Z}\frac{\sin^2 x_1}{(x_1-\pi n)^2}=\sum_{n\in\Z}\frac{\sin^2 y_2}{(y_2-\pi n)^2}=1.
$$
This implies the result. $\bl$

\medskip

{\bf Proof of Theorem \ref{SMdd}.} The result follows from Theorem \ref{fs},
because 
$$
\|\Phi\|_{\fM(E_1,E_2)}\le\|\Phi\|_{C_{\rm b}(\C)\hat\otimes_{\rm h}C_{\rm b}(\C)}
$$
for every $\Phi\in C_{\rm b}(\C)\hat\otimes_{\rm h}C_{\rm b}(\C)$ and for every Borel spectral
measures $E_1$ and $E_2$ on $\C$ (see \S\,\ref{dois}). $\bl$

\medskip

{\bf Remark.} The proof of Theorem \ref{SMdd} given above is based on the 
representation of \rf{haa1}. It is also possible to prove this theorem by using integral representation \rf{haa2} and estimate the norm in the space of Schur multipliers in terms of \rf{bs}.

\

\section{\bf Operator Lipschitzness and preservation of operator ideals}
\setcounter{equation}{0}
\label{OLpoi}

\

In this section we show that functions in the Besov space $B_{\be1}^1\big(\R^2\big)$
are operator Lipschitz. We also show that if $f\in B_{\be1}^1\big(\R^2\big)$, then 
$$
N_1-N_2\in\fI\quad\Longrightarrow\quad f(N_1)-f(N_2)\in\fI,
$$ 
whenever $\fI$ is a quasinormed operator ideal with majorization property.
In particular, this is true if $\fI=\bS_1$. 

Recall that in the case $\fI=\bS_1$ one cannot replace the Besov class $B_{\be1}^1(\R^2)$ with the Lipschitz class. Indeed, even in the case of self-adjoint operators a Lipschitz function $f$ on $\R$ does not possess the property
$$
A-B\in\bS_1\quad\Longrightarrow\quad f(A)-f(B)\in\bS_1.
$$
This was shown for the first time in \cite{F2}. Later necessary conditions were found in \cite{Pe2} and \cite{Pe3} that also show that Lipschitzness is not sufficient. 

The following lemma is an immediate consequence Theorems \ref{SMdd} and \ref{doin}.

\begin{lem}
\label{oLs}
Let $f$ be a function in $C_{\rm b}(\R^2)$ such that 
$$
\supp\F f\subset\{\z\in\C:~|\z|\le\s\},\quad\s>0.
$$
If $N_1$ and $N_2$ are normal operators, then 
$$
\|f(N_1)-f(N_2)\|\le\const\s\|f\|_{L^\be}\|N_1-N_2\|.
$$
\end{lem}

\begin{thm}
\label{oL}
Let $f$ belong to the Besov space $B_{\be1}^1\big(\R^2\big)$
and let $N_1$ and $N_2$ be normal operators whose difference is a bounded operator. Then {\em\rf{if}} holds and
$$
\|f(N_1)-f(N_2)\|\le\const\|f\|_{B_{\be1}^1(\R^2)}\|N_1-N_2\|.
$$
\end{thm}

\Pf It follows from Lemma \ref{oLs} that
\begin{align*}
\|f(N_1)-f(N_2)\|&\le\sum_{n\in\Z}\big\|f_n(N_1)-f_n(N_2)\big\|\\[.2cm]
&\le\const\sum_{n\in\Z}2^n\|f_n\|_{L^\be}\|N_1-N_2\|
\le\const\|f\|_{B_{\be1}^1(\R^2)}\|N_1-N_2\|
\end{align*} 
(see the definition of $B_{\be1}^1\big(\R^2\big)$ in \S\,\ref{Fsp}). $\bl$

In other words, functions in $B_{\be1}^1\big(\R^2\big)$ must be operator Lipschitz.

We can obtain similar results for operator ideals.

\begin{lem}
\label{oLis}
Let $\fI$ be a quasinormed ideal of operators on Hilbert space that has majorization property and let $f$ be a function in $C_{\rm b}\big(\R^2\big)$ such that 
$$
\supp\F f\subset\{\z\in\C:~|\z|\le\s\},\quad\s>0.
$$
If $N_1$ and $N_2$ are normal operators such that $N_1-N_2\in\fI$, then 
$$
f(N_1)-f(N_2)\in\fI\quad\mbox{and}\quad
\|f(N_1)-f(N_2)\|_\fI\le c\,\s\|f\|_{L^\be}\|N_1-N_2\|_\fI
$$
for a numerical constant $c$.
\end{lem}

\Pf The result follows from Theorem \ref{SMdd} and from \rf{intn}. $\bl$

\begin{thm}
\label{yad}
Let $\fI$ be a quasinormed ideal of operators on Hilbert space that has majorization property and let $f$ belong to the Besov space $B_{\be1}^1\big(\R^2\big)$.
If $N_1$ and $N_2$ are normal operators such that $N_1-N_2\in\fI$.
Then $f(N_1)-f(N_2)\in\fI$ and
$$
\|f(N_1)-f(N_2)\|_\fI\le c\,\|f\|_{B_{\be1}^1(\R^2)}\|N_1-N_2\|_\fI
$$
for a numerical constant $c$.
\end{thm}

\Pf In the case where $\fI$ is a normed ideal the result is an immediate 
consequence of Lemma \ref{oLis}.
In particular, Theorem \ref{yad} is true for $\fI=\bS_1^l$. To complete the 
proof in the general
case it suffices to use the majorization property. $\bl$

\begin{cor}
\label{S1B}
There exists a positive number $c$ such that if
 $f\in B_{\be1}^1\big(\R^2\big)$ and let $N_1$ and $N_2$ are normal operators such that $N_1-N_2\in\bS_1$, 
then $f(N_1)-f(N_2)\in\bS_1$ and
$$
\|f(N_1)-f(N_2)\|_{\bS_1}\le c\,\|f\|_{B_{\be1}^1(\R^2)}\|N_1-N_2\|_{\bS_1}.
$$
\end{cor}

\

\section{\bf Operator H\"older functions and arbitrary moduli of continuity}
\setcounter{equation}{0}
\label{oHfamc}

\

Recall that $\a\in(0,1)$, the class $\L_\a\big(\R^2\big)$ of H\"older functions of order $\a$ is defined by:
$$
\L_\a\big(\R^2\big)\df\left\{f:~\|f\|_{\L_\a(\R^2)}=
\sup_{z_1\ne z_2}\frac{|f(z_1)-f(z_2)|}{|z_1-z_2|^\a}<\be\right\}.
$$
In this section we show that in contrast with the class of Lipschitz functions, a H\"older function of order $\a\in(0,1)$ must be {\it operator H\"older of order $\a$}.

We also consider in this section the more general case of functions in the space $\L_\o\big(\R^2\big)$, where $\o$ is an arbitrary modulus of continuity.

\begin{thm}
\label{oH}
There exists a positive number $c$ such that for every $\a\in(0,1)$ and every
$f\in\L_\a\big(\R^2\big)$,
\bay
\label{alfa}
\|f(N_1)-f(N_2)\|\le c\,(1-\a)^{-1}\|f\|_{\L_\a(\R^2)}\|N_1-N_2\|^\a.
\ey
for arbitrary normal operators $N_1$ and $N_2$.
\end{thm}

\Pf The proof is almost the same as the proof of Theorem 4.1 of \cite{AP2}
(see also Remark 1 following Theorem 4.1 in \cite{AP2}) for self-adjoint operators. All we need is the following:
\bay
\label{in1}
\big\|f_n(N_1)-f_n(N_2)\big\|\le\const2^n\|f_n\|_{L^\be}\|N_1-N_2\|,\quad n\in\Z,
\ey
and
\bay
\label{in2}
\|f_n\|_{L^\be}\le\const2^{-n\a}\|f\|_{\L_\a(\R^2)},\quad n\in\Z,
\ey
where the functions $f_n$ are defined by \rf{fn}.
We remind that \rf{in1} is a consequence of Lemma \ref{oLs}, while \rf{in2} is a special case of Theorem \ref{Vn}.

The deduction of inequality \rf{alfa} from \rf{in1} and \rf{in2} is exactly the same as in the proof of Theorem 4.1 of \cite{AP2}, in which inequality \rf{alfa} for self-adjoint operators is deduced from the corresponding analogs of inequalities \rf{in1} and \rf{in2}. $\bl$

Consider now more general classes of functions. Let $\o$ be a modulus of continuity. Recall that the class $\L_\o\big(\R^2\big)$ is defined by
$$
\L_\o\big(\R^2\big)\df\left\{f:~\|f\|_{\L_\o(\R^2)}=
\sup_{z_1\ne z_2}\frac{|f(z_1)-f(z_2)|}{\o(|z_1-z_2|)}<\be\right\}.
$$

As in the case of functions of one variable (see \cite{AP1}, \cite{AP2}), we define the function $\o_*$ by
\bay
\label{o*}
\o_*(x)\df x\int_x^\be\frac{\o(t)}{t^2}\,dt,\quad x>0.
\ey

\begin{thm}
\label{Lo}
There exists a positive number $c$ such that for every
modulus of continuity $\o$ and every $f\in\L_\o\big(\R^2\big)$,
\bay
\label{omega}
\|f(N_1)-f(N_2)\|\le c\,\|f\|_{\L_\o(\R^2)}\,\o_*\big(\|N_1-N_2\|\big)
\ey
for arbitrary normal operators $N_1$ and $N_2$.
\end{thm}

\Pf To prove Theorem \ref{Lo}, we need inequalities \rf{in1} and Theorem \ref{Vn}. 
The deduction of inequality \rf{omega} from \rf{in1} and Theorem \ref{Vn} is exactly the same as it was done in the proof of Theorem 7.1 of \cite{AP2} in the case of self-adjoint operators. $\bl$

\begin{cor}
\label{sle}
Let $\o$ be a modulus of continuity such that
$$
\o_*(x)\le\const\o(x),\quad x>0,
$$
and let $f\in\L_\o(\R^2)$. Then
$$
\|f(N_1)-f(N_2)\|\le\const\|f\|_{\L_\o(\R^2)}\,\o\big(\|N_1-N_2\|\big)
$$
for arbitrary normal operators $N_1$ and $N_2$.
\end{cor}

Theorem \ref{Lo} allows us to estimate $\|f(N_1)-f(N_2)\|$ for Lipschitz functions $f$ and normal operators $N_1$ and $N_2$ whose spectra are contained in a given compact convex subset of $\C$.

For a Lipschitz function $f$ on a subset $K$ of $\C$, the Lipschitz constant is, by definition,
$$
\|f\|_{\Li}\df\sup\left\{\frac{|f(\z_1)-f(\z_2)|}{|\z_1-\z_2|}:~\z_1,\,\z_2\in K,~\z_1\ne\z_2\right\}.
$$
For a Lipschitz function $f$ on a compact convex subset $K$ of $\C$, we extend it to $\C$ by the formula 
\bay
\label{ras}
f(\z)\df f(\z_\sharp), 
\ey
where $\z_\sharp$ is the closest point to $\z$ in $K$. It is easy to see that
the Lipschitz constant of this extension does not change.

\begin{thm}
\label{ogra}
Let $N_1$ and $N_2$ be normal operators whose spectra are contained in 
a compact convex set $K$ and let $f$ be a Lipschitz function on $K$.
Then
\bay
\label{lipots}
\big\|f(N_1)-f(N_2)\big\|\le
\const\|f\|_{\Li}\|N_1-N_2\|\left(1+\log\frac{d}{\|N_1-N_2\|}\right),
\ey
where $d$ is the diameter of $K$.
\end{thm}

\Pf Without loss of generality, we may assume that $\|f\|_\Li=1$. Let us extend $f$ to $\C$ by formula \rf{ras}. Define the modulus of continuity $\o$ by
$$
\o(\d)=\left\{\begin{array}{ll}\d,&\d\le d,\\[.2cm]d,&\d>d.
\end{array}\right.
$$
Clearly, $f\in\L_\o\big(\R^2\big)$ and $\|f\|_{\L_\o(\R^2)}\le\|f\|_\Li$.
We have
$$
\o_*(\d)=\d\int_\d^d\frac{dt}t+\d d\int_d^\be\frac{dt}{t^2}
=\d\log\frac{d}{\d}+\d,\quad\d\le d,
$$
where $\o_*$ is defined by \rf{o*}. Now inequality \rf{lipots} follows immediately from Theorem \ref{Lo}. $\bl$

\

\section{\bf Perturbations of class $\bS_p$ and more general operator ideals}
\setcounter{equation}{0}
\label{pocS}

\

In this section we obtain sharp estimates for $f(N_1)-f(N_2)$ in the case when \lb$f\in\L_\a\big(\R^2\big)$, $0<\a<1$, and $N_1$ and $N_2$ are normal operators such whose difference belong to Schatten--von Neumann classes $\bS_p$. We also obtain more general results in the case when the difference of the normal operators belongs to operator ideals.

Let us first state the result for Schatten--von Neumann classes.

\begin{thm}
\label{sp}
Let $0<\a<1$ and $1<p<\be$. Then there exists a positive number $c$ such that
for every $f\in\L_\a\big(\R^2\big)$ and for arbitrary normal operators $N_1$ and $N_2$ with $N_1-N_2\in\bS_p$, the operator $f(N_1)-f(N_2)$ belongs to $\bS_{p/\a}$ and
the following inequality holds:
$$
\big\|f(N_1)-f(N_2)\big\|_{\bS_{p/\a}}\le c\,\|f\|_{\L_\a(\R)}\|N_1-N_2\|^\a_{\bS_p}.
$$
\end{thm}

We discuss the case $p=1$ after the proof of Theorem \ref{S1}.

Theorem \ref{sp} is an immediate consequence of a more general result for operator ideals, see Theorem \ref{pere} below.

To proceed to operator ideals, we start with the ideals $\bS_p^l$.
Recall that for $l\ge0$ and $p\ge1$, the normed ideal $\bS_p^l$ consists of all bounded linear operators equipped with the norm
$$
\|T\|_{\bS_p^l}\df\left(\sum_{j=0}^l\big(s_j(T)\big)^p\right)^{1/p}.
$$

\begin{thm}
\label{gl}
Let $0<\a<1$. Then there exists a positive number $c>0$ such that for every
$l\ge0$,  $p\in[1,\be)$,  $f\in\L_\a\big(\R^2\big)$, and for arbitrary normal operators $N_1$ and $N_2$ on Hilbert space with bounded $N_1-N_2$, the following inequality holds:
$$
s_j\big(f(N_1)-f(N_2)\big)\le c\,\|f\|_{\L_\a(\R^2)}(1+j)^{-\a/p}\|N_1-N_2\|_{\bS_p^l}^\a
$$
for every $j\le l$.
\end{thm}

\Pf The proof is almost the same as the proof of Theorem 5.1 of \cite{AP3}.
To be able to apply the reasonings given in the proof of Theorem 5.1 of \cite{AP3}, we need inequality \rf{in2} and the following inequality:
\bay
\label{inpl}
\big\|f_n(N_1)-f_n(N_2)\big\|_{\bS_p^l}\le\const2^n\|f_n\|_{L^\be}\|N_1-N_2\|_{\bS_p^l},\quad n\in\Z,
\ey
where the functions $f_n$ are defined by \rf{fn}. Inequality \rf{inpl} is an immediate consequence of Lemma \ref{oLis}. All the details can be found in the proof of Theorem 5.1 of \cite{AP3}. $\bl$ 

\begin{thm}
\label{S1}
Let $0<\a<1$. Then there exists a positive number $c>0$ such that for every
$f\in\L_\a\big(\R^2\big)$ and arbitrary normal operators $N_1$ and $N_2$ on Hilbert space with $N_1-N_2\in\bS_1$, the operator $f(N_1)-f(N_2)$ belongs to 
$\bS_{\frac1\a,\be}$ and the following inequality holds:
$$
\big\|f(N_1)-f(N_2)\big\|_{\bS_{\frac1\a,\be}}\le c\,\|f\|_{\L_\a(\R^2)}
\|N_1-N_2\|_{\bS_1}^\a.
$$
\end{thm}

\Pf As in the case of self-adjoint operators (see Theorem 5.2 of \cite{AP3}), this is an immediate consequence of Theorem \ref{gl} in the case $p=1$. $\bl$

Note that the assumptions of Theorem \ref{S1} do not imply that 
$f(N_1)-f(N_2)\in\bS_{1/\a}$. This is not true even in the case when $N_1$ and $N_2$ are self-adjoint operators. This was proved in \cite{AP3}. Moreover, in \cite{AP3} a necessary condition on the function $f$ on $\R$ was found
for
$$
f(A)-f(B)\in\bS_{1/\a},\quad\mbox{whenever}\quad
A=A^*,~B=B^*\qm A-B\in\bS_1.
$$
That necessary condition is based on the $\bS_p$ criterion for Hankel operators (\cite{Pe1} and \cite{Pe4}, Ch. 6) and shows that the condition
$f\in\L^\a(\R)$ is not sufficient.

The following result ensures that the assumption that $N_1-N_2\in\bS_1$
for normal operators $N_1$ and $N_2$
implies that $f(N_1)-f(N_2)\in\bS_{1/\a}$ under a slightly more restrictive assumption on $f$.

\begin{thm}
\label{S1s}
Let $0<\a\le1$. Then there exists a positive number $c>0$ such that for every
$f\in B_{\be1}^\a\big(\R^2\big)$ and arbitrary normal operators $N_1$ and $N_2$ on Hilbert space with $N_1-N_2\in\bS_1$, the operator $f(N_1)-f(N_2)$ belongs to 
$\bS_{1/\a}$ and the following inequality holds:
$$
\big\|f(N_1)-f(N_2)\big\|_{\bS_{1/\a}}\le c\,\|f\|_{B_{\be1}^\a(\R^2)}
\|N_1-N_2\|_{\bS_1}^\a.
$$
\end{thm}

Note that in the case $\a=1$ turns into Corollary \ref{S1B}.

\medskip

{\bf Proof of Theorem \ref{S1s}.} Again, if we apply Lemma \ref{oLis}, the proof is practically the same as the proof of Theorem 5.3 in \cite{AP3}. $\bl$

\begin{thm}
\label{sigma}
Let $0<\a<1$. Then there exists a positive number $c>0$ such that for every
$f\in\L_\a\big(\R^2\big)$ and arbitrary normal operators $N_1$ and $N_2$ on Hilbert space with bounded $N_1-N_2$, the following inequality holds:
$$
s_j\Big(\big|f(N_1)-f(N_2)\big|^{1/\a}\Big)\le c\,\|f\|_{\L_\a(\R)}^{1/\a}
\s_j(N_1-N_2),\quad j\ge0.
$$
\end{thm}

Recall that the numbers $\s_j(N_1-N_2)$
defined by \rf{sin}.

\medskip

\Pf As in the case of self-adjoint operators (see \cite{AP3}), it suffices to apply Theorem \ref{gl} with $l=j$ and $p=1$. $\bl$

Now we are in a position to obtain a general result in the case 
$f\in\L_\a\big(\R^2\big)$ and $N_1-N_2\in\fI$ for an arbitrary quasinormed ideal $\fI$ with upper Boyd index less than 1. Recall that the number $\bs{C}_\fI$ is defined in \S\,\ref{Opid}.

\begin{thm}
\label{osn}
Let $0<\a<1$. Then there exists a positive number $c>0$ such that for every
$f\in\L_\a\big(\R^2\big)$, for an arbitrary quasinormed ideal $\fI$ with $\b_\fI<1$, and for arbitrary normal operators $N_1$ and $N_2$ on Hilbert space with $N_1-N_2\in\fI$,
the operator $\big|f(N_1)-f(N_2)\big|^{1/\a}$ belongs to $\fI$
and the following inequality holds:
$$
\Big\|\big|f(N_1)-f(N_2)\big|^{1/\a}\Big\|_\fI\le c\,\bs{C}_\fI
\|f\|_{\L_\a(\R^2)}^{1/\a}\|N_1-N_2\|_\fI.
$$
\end{thm}

\Pf The proof is almost the same as the proof of Theorem 5.5 in \cite{AP3}. $\bl$

We can reformulate Theorem \ref{osn} in the following way.

\begin{thm}
\label{pere}
Under the hypothesis of Theorem {\em\ref{osn}}, the operator $f(N_1)-f(N_2)$
belongs to $\fI^{\{1/\a\}}$ and
$$
\big\|f(N_1)-f(N_2)\big\|_{\fI^{\{1/\a\}}}\le c^\a\,\bs{C}_\fI^\a\|f\|_{\L_\a(\R^2)}
\|N_1-N_2\|^\a_\fI.
$$
\end{thm}

The following result is a consequence of Theorem \ref{osn}.

\begin{thm}
\label{spl}
Let $0<\a<1$ and $1<p<\be$. Then there exists a positive number $c$ such that
for every $f\in\L_\a\big(\R^2\big)$, every $l\in\Z_+$, and arbitrary normal operators $N_1$ and $N_2$ with bounded $N_1-N_2$, the following inequality holds:
$$
\sum_{j=0}^l\left(s_j\Big(\big|f(N_1)-f(N_2)\big|^{1/\a}\Big)\right)^p\le
c\,\|f\|_{\L_\a(\R^2)}^{p/\a}\sum_{j=0}^l\big(s_j(N_1-N_2)\big)^p.
$$
\end{thm}

\Pf As in the case of self-adjoint operators (see \cite{AP3}), the result immediately follows from Theorem \ref{osn} from \rf{ClI}.
$\bl$

\

\section{\bf Commutators and quasicommutators}
\setcounter{equation}{0}
\label{Caq}

\

In this section we obtain estimates for quasicommutators
$f(N_1)R-Rf(N_2)$, where $N_1$ and $N_2$ are normal operators and 
$R$ is a bounded linear operator. In the special case when $R=I$ we arrive at the problem of estimating $f(N_1)-f(N_2)$ that we have discussed in previous sections. On the other hand, in the special case when $N_1=N_2$ we have the problem of estimating commutators $f(N)R-Rf(N)$.

It turns out, however, that it is impossible to obtain estimates of
$\|f(N_1)R-Rf(N_2)\|$ in terms of $\|N_1R-RN_2\|$. This cannot be done even for the function $f(z)=\bar z$.

Though the well-known Fuglede-Putnam theorem says that the equality
$N_1R=RN_2$ for a bounded operator $R$ and normal operators $N_1$ and $N_2$ implies that $N_1^*R=RN_2^*$, the smallness of $N_1R-RN_2$ does not imply
the smallness of $N_1^*R-RN_2^*$. 

Indeed, it follows from Corollary 4.3 of \cite{JW} that for every $\e>0$ there exists a bounded normal operator $N$ and operator $R$ of norm 1 such that
$$
\|NR-RN\|<\e\quad\mbox{but}\quad\|N^*R-RN^*\|\ge1.
$$

The results of \cite{JW} also imply that
if $f\in C(\C)$ and 
$$
\|f(N)Q-Qf(N)\|\le\const\|NQ-QN\|
$$
for all bounded operators $Q$ and
bounded normal operators $N$, then $f$ is a linear function, i.e., $f(z)=az+b$ for some $a,\,b\in\C$.

{\it In this section we obtain estimates for quasicommutators $f(N_1)R-Rf(N_2)$
in terms of the quasicommutators $N_1R-RN_2$ and $N_1^*R-RN_2^*$.}

%Indeed, the results of \cite{JW} implies that $f$ is differentiable as a 
%function
%of a complex variable everywhere. Thus, $f$ is an entire function. It 
%remains to
%note that $|f(z)-f(w)|\le\const|z-w|$ for all $z, w\in\C$.

Let us explain what we mean by the boundedness of $N_1R-RN_2$
for not necessarily bounded normal operators $N_1$ and $N_2$.

We say that the {\it operator $N_1R-RN_2$ is bounded} if
$R(\cd_{N_2})\subset \cd_{N_1}$ and 
$$
\|N_1Ru-RN_2u\|\le\const\|u\|\quad\mbox{for every}\quad u\in \cd_{N_2}.
$$
Then there exists a unique bounded operator $K$ such that
$Ku=N_1Ru-RN_2u$ for all $u\in \cd_{N_2}$. In this case we write $K=N_1R-RN_2$. Thus $N_1R-RN_2$ is bounded if and only if
\bay
\label{MN}
\big|(Ru,N_1^*v)-(N_2u,R^*v)\big|\le \const\|u\|\cdot\|v\|
\ey
for every $u\in \cd_{N_2}$ and $v\in \cd_{N_1^*}=\cd_{N_1}$. It is easy to see that
$N_1R-RN_2$ is bounded if and only if $N_2^*R^*-R^*N_1^*$ is bounded,
and $(N_1R-RN_2)^*=-(N_2^*R^*-R^*N_1^*)$.
In particular, we write $N_1R=RN_2$ if $R(\cd_{N_2})\subset \cd_{N_1}$ and $N_1Ru=RN_2u$ for every $u\in \cd_{N_2}$.
We say that $\|N_1R-RN_2\|=\be$ if $N_1R-RN_2$ is not a bounded operator.

We need the following observation:

\medskip

{\bf Remark.} Suppose that $N_1^*$ is the closure of an operator
$N_\flat$ and $N_2$ is the closure of an operator $N_\sharp$. Suppose that inequality
\rf{MN} holds for all $u\in \cd_{N_\sharp}$ and $v\in\cd_{N_\flat}$.
Then it holds for all $u\in \cd_{N_2}$ and $v\in \cd_{N_1}$.

\begin{thm}
\label{cqc}
Let $f$ be a function in $C_{\rm b}\big(\R^2)$ whose Fourier transform $\F f$ has compact support. Suppose that $R$ is a bounded linear operator, $N_1$ and $N_2$ are normal operators such that the operators $N_1R-RN_2$ and
$N^*_1R-RN^*_2$ are bounded. Then
\begin{align}
\label{doic}
f(N_1)R-Rf(N_2)&=\iint\limits_{\C^2}\big(\dg_yf\big)(z_1,z_2)\,
dE_1(z_1)(B_1R-RB_2)\,dE_2(z_2)\nonumber\\[.2cm]
&+\iint\limits_{\C^2}\big(\dg_xf\big)(z_1,z_2)\,
dE_1(z_1)(A_1R-RA_2)\,dE_2(z_2).
\end{align}
\end{thm}

\Pf The proof is similar to the proof of Theorem \ref{doin},
Consider first the case when $N_1$ and $N_2$ are bounded operators. Put 
$$
d=\max\big\{\|N_1\|,\|N_2\|\big\}\quad\mbox{and}\quad D\df\{\z\in\C:~|\z|\le d\}.
$$

By Theorem \ref{SMdd}, both $\dg_yf$ and $\dg_xf$ are Schur multipliers.
We have
\begin{align*}
\iint\limits_{\C^2}&\big(\dg_yf\big)(z_1,z_2)\,
dE_1(z_1)(B_1R-RB_2)\,dE_2(z_2)\\[.2cm]
&=\iint\limits_{D\times D}\big(\dg_yf\big)(z_1,z_2)\,
dE_1(z_1)(B_1R-RB_2)\,dE_2(z_2)\\[.2cm]
&=\iint\limits_{D\times D}\big(\dg_yf\big)(z_1,z_2)\,dE_1(z_1)B_1R\,dE_2(z_2)
-\iint\limits_{D\times D}\big(\dg_yf\big)(z_1,z_2)\,dE_1(z_1)RB_2\,dE_2(z_2)
\\[.2cm]
&=\iint\limits_{D\times D}y_1\big(\dg_yf\big)(z_1,z_2)\,dE_1(z_1)R\,dE_2(z_2)
-\iint\limits_{D\times D}y_2\big(\dg_yf\big)(z_1,z_2)\,dE_1(z_1)R\,dE_2(z_2)
\\[.2cm]
&=\iint\limits_{D\times D}(y_1-y_2)\big(\dg_yf\big)(z_1,z_2)\,dE_1(z_1)R\,dE_2(z_2)\\[.2cm]
&=\iint\limits_{D\times D}\big(f(x_1,y_1)-f(x_1,y_2)\big)\,dE_1(z_1)R\,dE_2(z_2).
\end{align*}
Similarly,
\begin{align*}
\iint\limits_{\C^2}&\big(\dg_xf\big)(z_1,z_2)\,
dE_1(z_1)(A_1R-RA_2)\,dE_2(z_2)\\[.2cm]
&=\iint\limits_{D\times D}\big(f(x_1,y_2)-f(x_2,y_2)\big)\,dE_1(z_1)R\,dE_2(z_2).
\end{align*}
It follows that
\begin{align*}
\iint\limits_{\C^2}&\big(\dg_yf\big)(z_1,z_2)\,
dE_1(z_1)(B_1R-RB_2)\,dE_2(z_2)\\[.2cm]
+&\iint\limits_{\C^2}\big(\dg_xf\big)(z_1,z_2)\,
dE_1(z_1)(A_1R-RA_2)\,dE_2(z_2)\\[.2cm]
&=\iint\limits_{D\times D}\big(f(x_1,y_1)-f(x_2,y_2)\big)\,dE_1(z_1)R\,dE_2(z_2)\\[.2cm]
&=\iint\limits_{D\times D}f(x_1,y_1)\,dE_1(z_1)R\,dE_2(z_2)
-\iint\limits_{D\times D}f(x_2,y_2)\,dE_1(z_1)R\,dE_2(z_2)\\[.2cm]
&=f(N_1)R-Rf(N_2).
\end{align*} 

In the general case we use the same approximation procedure as in the proof of Theorem \ref{doin}. $\bl$

As in the case of differences $f(N_1)-f(N_2)$, we can extend Theorem \ref{cqc} to functions $f$ in $B_{\be1}^1\big(\R^2\big)$.

\begin{thm}
\label{dlyac}
Let $N_1$ and $N_2$ be normal operators and let $R$ be a bounded linear operator such that the quasicommutators $N_1R-RN_2$ and
$N^*_1R-RN^*_2$ are bounded.
Then {\em\rf{doic}} holds for every $f\in B_{\be1}^1\big(\R^2\big)$.
\end{thm}

\Pf The proof is almost the same as the proof of Theorem \ref{dlya}. $\bl$

Theorem \ref{dlyac} allows us to generalize all the results of Sections
\ref{OLpoi}, \ref{oHfamc}, and \ref{pocS} to the case of quasicommutators. 
We state some of them. The proofs of the theorems stated below is exactly the same as the proofs of the corresponding results in Sections
\ref{OLpoi}--\ref{pocS}.

\begin{thm}
\label{Besov}
There exists a positive number $c$ such that for every normal operators
$N_1$ and $N_2$, every bounded linear operator $R$ and an arbitrary function $f$ in $B_{\be1}^1\big(\R^2\big)$
the following inequality holds:
$$
\|f(N_1)R-Rf(N_2)\|\le c\|f\|_{B_{\be1}^1(\R^2)}
\max\big\{\|N_1R-RN_2\|,\|N^*_1R-RN^*_2\|\big\}.
$$
\end{thm}

\begin{thm}
\label{Holder}
Let $0<\a<1$. Then there exists $c>0$ such that for every $f\in\L_\a\big(\R^2\big)$, for arbitrary normal operators $N_1$ and $N_2$ and a bounded operator $R$ the following inequality holds:
$$
\|f(N_1)R-Rf(N_2)\|\le c\,\|f\|_{\L_\a(\R^2)}
\max\big\{\|N_1R-RN_2\|,\|N^*_1R-RN^*_2\|\big\} ^\a\|R\|^{1-\a}.
$$
\end{thm}

\begin{thm}
\label{modn}
There exists $c>0$ such that for every
modulus of continuity $\o$, for every $f\in\L_\o\big(\R^2\big)$, for arbitrary normal operators $N_1$ and $N_2$, and a bounded nonzero
operator $R$ the following inequality holds:
$$
\|f(N_1)R-Rf(N_2)\|\le c\,\|f\|_{\L_\o(\R^2)}
\|R\|\,\,\o_*\!
\left(\frac{\max\big\{\|N_1R-RN_2\|,\|N^*_1R-RN^*_2\|\big\}}{\|R\|}\right).
$$
\end{thm}

The next result shows that in the case $N_1R-RN_2\in\bS_p$, $1<p<\be$, and \lb$f\in\L_\a\big(\R^2\big)$, $0<\a<1$, we can estimate $\big\|f(N_1)R-Rf(N_2)\big\|_{\bS_{p/\a}}$ in terms of \lb$\|N_1R-RN_2\|_{\bS_p}$, we do not need $\|N^*_1R-RN^*_2\|_{\bS_p}$.

\begin{thm}
\label{spqc}
Let $0<\a<1$ and $1<p<\be$. Then there exists a positive number $c$ such that
for every $f\in\L_\a\big(\R^2\big)$, for arbitrary normal operators $N_1$ and $N_2$ and a bounded operator $R$ with $N_1R-RN_2\in\bS_p$ and $N^*_1R-RN^*_2\in\bS_p$, the operator $f(N_1)R-Rf(N_2)$ belongs to $\bS_{p/\a}$ and
the following inequality holds:
$$
\big\|f(N_1)R-Rf(N_2)\big\|_{\bS_{p/\a}}\le c\,\|f\|_{\L_\a(\R^2)}
\|N_1R-RN_2\|_{\bS_p}^\a.
$$
\end{thm}

\Pf In the same way as in the proof of Theorem \ref{sp}, we can prove that
$$
\big\|f(N_1)R-Rf(N_2)\big\|_{\bS_{p/\a}}\le c\,\|f\|_{\L_\a(\R^2)}
\max\left\{\|N_1R-RN_2\|_{\bS_p},\|N^*_1R-RN^*_2\|_{\bS_p}\right\}^\a.
$$
The result follows from the well-known inequality:
\bay
\label{1pbe}
\|N^*_1R-RN^*_2\|_{\bS_p}\le\const \|N_1R-RN_2\|_{\bS_p},\quad1<p<\be,
\ey
see \cite{AD} and \cite{S}. $\bl$

Note that inequality \rf{1pbe} does not hold for $p=1$, see \cite{KS}. Thus to obtain analogs of Theorems \ref{S1} and \ref{S1s}, we have to estimate the quasicommutators
$f(N_1)R-Rf(N_2)$ in terms of both $N_1R-RN_2$ and $N^*_1R-RN_2^*$. Let us state e.g., the analog of Theorem \ref{S1s}.

\begin{thm}
\label{B1aqc}
Let $0<\a<1$. Then there exists a positive number $c$ such that
for every $f\in B_{\be1}^\a\big(\R^2\big)$, for arbitrary normal operators $N_1$ and $N_2$ and a bounded operator $R$ with $N_1R-RN_2\in\bS_1$ and $N^*_1R-RN^*_2\in\bS_1$, the operator $f(N_1)R-Rf(N_2)$ belongs to $\bS_{1/\a}$ and
the following inequality holds:
$$
\big\|f(N_1)R-Rf(N_2)\big\|_{\bS_{1/\a}}\le c\,\|f\|_{B_{\be1}^\a(\R^2)}
\max\left\{\|N_1R-RN_2\|_{\bS_1},\|N^*_1R-RN^*_2\|_{\bS_1}\right\}^\a.
$$
\end{thm}

The proof is almost the same as the proof of Theorem \ref{S1s}.

\

\

\footnotesize
\noindent
\begin{tabular}{p{4.7cm}p{4cm}p{5cm}}
A.B. Aleksandrov & V.V. Peller & D.S. Potapov \& F.A. Sukochev\\
St.Petersburg Branch & Department of Mathematics & School of Mathematics \& Statistics\\
Steklov Institute of Mathematics  & Michigan State University &University of NSW\\
Fontanka 27  & East Lansing & Kensington NSW 2052\\
191023 St-Petersburg &Michigan 48824&Australia \\
Russian Federation & USA
\end{tabular}


\begin{thebibliography}{99}
\label{bibl}
\bibitem[AD]{AD} {\sc A. Abdessemed and E. B. Davies}, {\em Some commutator estimates in the Schatten classes},
J. Lond. Math. Soc. {\bf 39} (1989), 299--308.

\bibitem[AP1]{AP1} {\sc A.B. Aleksandrov} and {\sc V.V. Peller}, {\em Functions  of perturbed operators},
 C.R. Acad. Sci. Paris, S\'er I {\bf347} (2009), 483--488.

\bibitem[AP2]{AP2}  {\sc A.B. Aleksandrov} and {\sc V.V. Peller},  {\em Operator H\"older--Zygmund functions}, Advances in Math.
{\bf224} (2010), 910-–966.

\bibitem[AP3]{AP3}  {\sc A.B. Aleksandrov} and {\sc V.V. Peller},  {\em Functions of operators under perturbations of
class $\bS_p$}, J. Funct. Anal. {\bf258} (2010), 3675--3724.

\bibitem[AP4]{AP4}  {\sc A.B. Aleksandrov} and {\sc V.V. Peller},  {\em Functions of perturbed unbounded self-adjoint operators.
Operator Bernstein type inequalities}, Indiana Univ. Math. J., in press.

\bibitem[AP5]{AP5}  {\sc A.B. Aleksandrov} and {\sc V.V. Peller},
{\em Functions of perturbed dissipative operators}, to appear.

\bibitem[APPS]{APPS} {\sc A.B. Aleksandrov, V.V. Peller, D. Potapov}, and
{\sc F. Sukochev}, {\em Functionsof perturbednormal operators}, C.R. Acad. Sci. Paris, S\'er I {\bf348} (2010), 553--558.


\bibitem[Be]{Be} {\sc G. Bennett}, {\em Schur multipliers}, Duke Math. J. {\bf44} (1977), 603--639.

\bibitem[BS1]{BS1} {\sc M.S. Birman} and {\sc M.Z. Solomyak},
{\em Double Stieltjes operator integrals},
Problems of Math. Phys., Leningrad. Univ. {\bf1} (1966), 33--67 (Russian).
English transl.: Topics Math. Physics {\bf1} (1967), 25--54, Consultants Bureau Plenum
Publishing Corporation, New York.

\bibitem[BS2]{BS2} {\sc M.S. Birman} and {\sc M.Z. Solomyak},
 {\em Double Stieltjes operator integrals. II},
 Problems of Math. Phys., Leningrad. Univ. {\bf2} (1967), 26--60 (Russian).
English transl.: Topics Math. Physics {\bf2} (1968), 19--46, Consultants Bureau Plenum
Publishing Corporation, New York.

\bibitem[BS3]{BS3} {\sc M.S. Birman} and {\sc M.Z. Solomyak},
{\em Double Stieltjes operator integrals. III},
Problems of Math. Phys., Leningrad. Univ. {\bf6} (1973), 27--53 (Russian).

\bibitem[BS4]{BS0} {\sc M.Sh. Birman and M.Z. Solomyak}, {\it Spectral theory of 
selfadjoint operators in Hilbert space},
Mathematics and its Applications (Soviet Series), 
D. Reidel Publishing Co., Dordrecht, 1987.

\bibitem[BS5]{BS} {\sc M.S. Birman} and {\sc M.Z. Solomyak},
{\em Tensor product of a finite number of spectral measures is always a spectral measure}, Integral Equations Operator Theory {\bf24} (1996),  179--187.

\bibitem[BS6]{BS4} {\sc M.S. Birman} and {\sc M.Z. Solomyak}, {\em Double operator integrals
in Hilbert space},  Int. Equat. Oper. Theory  {\bf47}  (2003), 131--168.

\bibitem[DK]{DK} {\sc Yu.L. Daletskii} and {\sc S.G. Krein}, {\it Integration and differentiation of
functions of Hermitian operators and application to the theory of perturbations} (Russian), Trudy Sem.
Functsion. Anal., Voronezh. Gos. Univ. {\bf1} (1956), 81--105.

\bibitem[F1]{F1}  {\sc Yu.B. Farforovskaya}, {\em  The connection of the Kantorovich-Rubinshtein metric for spectral resolutions of selfadjoint operators with functions of operators},
Vestnik Leningrad. Univ.  {\bf19}  (1968), 94--97. (Russian).

\bibitem[F2]{F2}  {\sc Yu.B. Farforovskaya}, {\em An example of a Lipschitzian function of selfadjoint
operators that yields a nonnuclear increase under a nuclear perturbation}.  Zap. Nauchn. Sem.
Leningrad. Otdel. Mat. Inst. Steklov. (LOMI)  {\bf30}  (1972), 146--153 (Russian).

%\bibitem[5]{F2}  {\sc Yu.B. Farforovskaya},
%{\em An estimate of the norm of $\mid f(B)-f(A)\mid $ for selfadjoint operators $A$ and $B$},
%Zap. Nauchn. Sem. Leningrad. Otdel. Mat. Inst.
%Steklov. (LOMI)  {\bf56}  (1976), 143--162 (Russian).

\bibitem[FN1]{FN1} {\sc Yu.B. Farforovskaya} and {\sc L.N. Nikolskaya},
{\em An inequality for commutators of normal operators}, Acta Sci. Math. (Szeged) {\bf71} (2005), 751--765.

\bibitem[FN2]{FN2} {\sc Yu.B. Farforovskaya} and {\sc L.N. Nikolskaya}, {\em Operator H\"olderness of
H\"older functions} (Russian), to appear in
 Algebra i Analiz.
 
\bibitem[GK]{GK} I.C. Gohberg and M.G. Krein, {\it Introduction to the theory
of linear nonselfadjoint operators in Hilbert space,} Nauka, Moscow, 1965.

English transl.: Amer. Math. Soc., Providence, RI, 1969.
 
\bibitem[JW]{JW} {\sc B.E. Johnson} and {\sc J.P. Williams}, {\em The range of a normal derivation}, Pacific J. Math. {\bf58} (1975), 105--122.

\bibitem[K]{K}  {\sc T. Kato}, {\em Continuity of the map $S\mapsto|S|$ for linear operators},
Proc. Japan Acad.  {\bf49}  (1973), 157--160.

\bibitem[KS]{KS} {\sc E. Kissin and V. S. Shulman}, {\em Classes of operator-smooth functions. III Stable functions and
Fuglede ideals}, Proc. Edinburgh Math. Soc. {\bf 48} (2005), 175--197.

\bibitem[L]{L}
{\sc B.~Ya. Levin}, {\em Lectures on entire functions},
Translation of Math. Monogr., vol. 150, 1996.

\bibitem[Mc]{Mc}
{\sc A. McIntosh}, {\em Counterexample to a question on commutators}, 
Proc. Amer. Math. Soc. {\bf29} (1971) 337--340.
 
\bibitem[Pee]{Pee} {\sc J. Peetre},
{\em New thoughts on Besov spaces}, Duke Univ. Press., Durham, NC, 1976.

\bibitem[Pe1]{Pe1} {\sc V.V.Peller}, {\em Hankel operators of class ${\bf S}_{p}$ 
and their applications (rational approximation, Gaussian processes, 
the problem of majorizing operators)}, Mat. Sbornik, 
{\bf 113} (1980), 538-581. 

English Transl. in Math. USSR Sbornik, {\bf 41}
(1982), 443-479.


\bibitem[Pe2]{Pe2} {\sc V.V. Peller},
{\em Hankel operators in the theory of perturbations of unitary and self-adjoint operators},
Funktsional. Anal. i Prilozhen. {\bf19:2}  (1985),
37--51 (Russian).
English transl.: Funct. Anal. Appl. {\bf19} (1985) , 111--123.

\bibitem[Pe3]{Pe3} {\sc V.V. Peller},
{\em Hankel operators in the perturbation theory of of unbounded self-adjoint operators}.
Analysis and partial differential equations,  529--544,
Lecture Notes in Pure and Appl. Math., {\bf122}, Dekker, New York, 1990.

\bibitem[Pe4]{Pe4} {\sc V.V. Peller}, {\em Hankel operators and their applications,}
Springer-Verlag, New York, 2003.

\bibitem[Pe5] {Pe5}  {\sc V.V. Peller}, {\em Differentiability of functions of contractions}, In: Linear and complex analysis, AMS Translations, Ser. 2 {\bf226}
(2009), 109--131, AMS, Providence.

\bibitem[Pi]{Pi} {\sc G. Pisier}, {\em Similarity problems and completely bounded maps},
Second, expanded edition. Includes the solution to ``The Halmos problem''. Lecture Notes in Mathematics, 1618. Springer-Verlag, Berlin, 2001.

\bibitem[S]{S} {\sc V.S. Shulman}, {\em Some remarks on the Fuglede-Weiss theorem,} Bull. London Math. Soc.
{\bf 28}(1996), 385--392.

%
%\bibitem[Pe6] {Pe6}  {\sc V.V. Peller}, {\em Differentiability of functions of contractions}, preprint, 2008.

%\bibitem[St]{S} {\sc V.V. Sten'kin}, {\em Multiple operator integrals},
%Izv. Vyssh. Uchebn. Zaved. Matematika {\bf4 (79)} (1977), 102--115 (Russian).

%English transl.: Soviet Math. (Iz. VUZ) {\bf21:4} (1977), 88--99.

%\bibitem[W]{W} {\sc H. Widom}, {\it When are differentiable functions differentiable?}, In: Linear and
%Complex Analysis Problem Book, Lect. Notes Math. {\bf 1043} (1984), 184--188.



%\bibitem[SNF]{SNF} {\sc B. Sz.-Nagy and C. Foias},
%{\it Harmonic analysis of operators on Hilbert
%space,} Akad\'{e}miai Kiad\'{o}, Budapest, 1970.


%\bibitem[SS]{SS} {\sc M.Z. Solomyak} and {\sc V.V. Sten'kin},
%{\em A certain class of multiple operator Stieltjes integrals},
%Problems of Math. Anal., No. 2: Linear Operators and Operator Equations (Russian), 122--134.
%Izdat. Leningrad. Univ., Leningrad, 1969.

\bibitem[T]{T} {\sc H. Triebel}, {\em Theory of function spaces,}
Monographs in Mathematics, {\bf78}, Birkh\"auser Verlag, Basel, 1983.


\end{thebibliography}
\end{document}